\documentclass[10pt,a4paper]{article}
  
\usepackage[T1]{fontenc}
\usepackage[utf8]{inputenc}
\usepackage{lmodern}
\usepackage[english]{babel}
\usepackage{amsmath,amsfonts,amssymb,amsthm}
\usepackage{pifont}
\usepackage{bm}
\usepackage{hyperref}
\usepackage{setspace}
\usepackage{thmtools}
\usepackage{stmaryrd}
\usepackage{fullpage}
\usepackage[shortlabels]{enumitem}
\usepackage[all]{xy}
\usepackage{xcolor}
\usepackage{tikz}

\usetikzlibrary{positioning} 

\definecolor{vi1}{RGB}{108,2,119}
\definecolor{v1}{RGB}{0,86,27}
\definecolor{link}{RGB}{175,25,25}
\hypersetup{colorlinks = true,allcolors=link}

\declaretheoremstyle[
  spaceabove=\topsep, spacebelow=\topsep,
  headfont=\bfseries\scshape,
  headpunct={.},
  notefont=\mdseries, notebraces={(}{)},
  bodyfont=\normalfont,
  postheadspace=1em,
  qed=\qedsymbol
]{mythmstyle}
\declaretheorem[style=mythmstyle,numbered=no,name=Proof]{dem}

\theoremstyle{plain}
\newtheorem{thm}{Theorem}[section]   
\newtheorem{lem}[thm]{Lemma}          
\newtheorem{prop}[thm]{Proposition}
\newtheorem{cor}[thm]{Corollary}

\theoremstyle{definition}
\newtheorem{defi}[thm]{Definition}

\newtheorem{rem}[thm]{Remark}
\newtheorem{rems}[thm]{Remarks}
\newtheorem{ex}[thm]{Example}
\newtheorem{exs}[thm]{Examples}
\newtheorem{nota}[thm]{Notation}

\title{Partial groups, pregroups and realisability of fusion systems}
\author{Nicolas Lemoine and R\'emi Molinier}

\numberwithin{equation}{subsection}
\numberwithin{figure}{section}

\newcommand{\Db}{\mathbb{D}}
\newcommand{\Fb}{\mathbb{F}}
\newcommand{\Nb}{\mathbb{N}}
\newcommand{\Wb}{\mathbb{W}}
\newcommand{\Zb}{\mathbb{Z}}

\newcommand{\Bc}{\mathcal{B}}
\newcommand{\Cc}{\mathcal{C}}
\newcommand{\Fc}{\mathcal{F}}
\newcommand{\F}{\mathcal{F}}
\newcommand{\Gc}{\mathcal{G}}
\newcommand{\Lc}{\mathcal{L}}
\newcommand{\Mc}{\mathcal{M}}
\newcommand{\Nc}{\mathcal{N}}
\newcommand{\Pc}{\mathcal{P}}
\newcommand{\Tc}{\mathcal{T}}

\newcommand{\ParG}{\mathbf{ParG}}
\newcommand{\PrG}{\mathbf{PrG}}
\newcommand{\Grps}{\mathbf{Grps}}
\newcommand{\sSet}{\mathbf{sSet}}
\newcommand{\Top}{\mathbf{Top}}
\newcommand{\triv}{{\mathbf{1}}}

\DeclareMathOperator{\Hom}{Hom}

\DeclareMathOperator{\Inj}{Inj}
\DeclareMathOperator{\Mor}{Mor}
\DeclareMathOperator{\Ob}{Obj}
\DeclareMathOperator{\Syl}{Syl}

\begin{document}
\maketitle

\begin{abstract}
In this article, we compare two different notions of partially defined group strutures, namely partial groups and pregroups, as introduced by Chermak and Stallings respectively. In particular we prove that the category of pregroups can be seen as a full subcategory of the category of partial groups. We also bring out some conjugation properties about elements and subgroups of finite order in pregroups and their universal groups. We then use these to investigate the question of realisability of fusion systems in finite pregroups.
\end{abstract}

\tableofcontents

\section*{Introduction}

Pregroups are generalisations of groups, where the product is only partially defined. Specifically, a \emph{pregroup} is a set $P$ equipped with a binary product $m\colon D\to P$, where $D$ is a subset of $P\times P$, subject to some group-like axioms (see Definition \ref{def:Pre} for details). Pregroups were introduced by Stallings \cite{St} in the 70's, as a tool to study amalgamated product of groups. A crucial property of these pregroups is that they can always be embedded into a group. More precisely, given a pregroup $P$ there is a group $U(P)$ containing $P$, called the \emph{universal group} of $P$, such that $P$ generates $U(P)$ in a strong way, implying that $U(P)$ has a solvable word problem when $P$ is finite. In an article published in 1987 \cite{Ri}, Rimlingler proved that under one finiteness condition on $P$, its universal group $U(P)$ is actually the fundamental group of a graph of groups whose edge and vertex groups are subgroups of $P$. He also proved that given a graph of groups of finite diameter, if all the edge maps are injective but not surjective, then its fondamental group is the universal group of a certain pregroup. Therefore, pregroups can be thought as combinatorial objects which encode graphs of groups. 

With other perspectives in mind, in an article published in 2013 \cite{Ch0}, Andrew Chermak introduced another generalisation of groups with a partially defined product : partial groups. He was interested in $p$-local structures of finite groups (i.e. how a group acts by conjugation on its $p$-subgroups) and more generally of fusion systems. A \emph{fusion system} over a finite $p$-group $S$ is a category whose objects are the subgroups of $S$ and whose morphism sets are formed by the conjugation maps induced by elements of $S$ together with other similar injective group homomorphisms (see Definition \ref{defF} for more details). The archetypal example of a fusion system is the one induced by a finite group acting by conjugation on one of its Sylow $p$-subgroups. If a fusion system satisfies a few more axioms ensuring that it ``behaves'' like one of these archetypal examples, we say that the fusion system is \emph{saturated}. Actually, general fusion systems give a very large class of mathematical objects, which strictly contains fusion systems induced by finite groups. There even exist saturated fusion systems, called \emph{exotic}, which cannot be obtained from a finite group $G$ containing $S$ as a Sylow $p$-subgroup. This leads to the question of the \emph{realisability} of saturated fusion systems: given a saturated fusion system over a finite $p$-group $S$, can it be obtained from a finite group containing $S$ as a Sylow $p$-subgroup? This question is also of interest (and leads to different answers) if we restate it in a larger framework: considering general fusion systems, dropping the Sylow $p$-subgroup hypothesis, or enlarging the context of realisability (possibly infinite groups, partial groups...).

A \emph{partial group}, as introduced by Chermak, is a set $\Mc$ together with a multivariate product defined on a subset of the set of words in $\Mc$ and which satisfies some axioms (see Definition \ref{def:PG} for details). One of Chermak's main achievements was to prove the existence and uniqueness of a so-called \emph{linking system} associated to a given saturated fusion system. In the way, he proved that there is a one-to-one correspondence between localities -- a certain type of partial groups that models $p$-local structures of finite groups -- and transporter systems (see Appendix A in  \cite{Ch0}).
 Linking systems and transporter systems are categories derived from a given saturated fusion system and they were introduced by Broto, Levi and Oliver in \cite{BLO2} and Oliver and Ventura in \cite{OV1} respectively. They used these objects to study saturated fusion systems, $p$-completed classifying spaces of finite groups and connections between them, by developing a theory of classifying spaces for saturated fusion systems. Localities give a more group-like point of view on these objects, which allows for instance the use of tools from group theory. 

Even though partial groups are combinatorial objects, they can be viewed as simplicial sets as highlighted by Broto and Gonzales \cite{BG}. For example, the geometric realisation of a locality has the homotopy type of the geometric realisation of the nerve of the associated transporter system. Therefore it opened another approach to study these spaces which are crucial in homotopy theory of fusion systems.

In this paper we establish connections between pregroups and partial groups, and we develop the question of the realisability of fusion systems in this generalised context. In Section \ref{Section1}, we give the basic definitions and properties concerning Chermak's partial groups, including the simplicial point of view. Section \ref{Section2} is dedicated to pregroups. After giving definitions, examples and some lemmas concerning conjugation properties and elements of finite order in pregroups, we prove that pregroups can be viewed as partial groups in a natural way. More precisely, there is a fully faithful functor from the category of pregroups $\PrG$ to the category of partial groups $\ParG$ (Proposition \ref{prop pregroup is partial group}). As a consequence, we get that the universal group $U(P)$ of a pregroup $P$ is isomorphic to the fundamental group of the geometric realisation $BP$ of the corresponding simplicial set $\mathcal{B}(P)$, called the classifying space of $P$, when $P$ is considered as a partial group. 

\begin{thm}[Corollary \ref{cor:universal=pi1} below] \label{thm:intro-universal=pi1}
Let $P$ be a pregroup. Then, $\pi_1(BP)\cong U(P)$.
\end{thm}

In Section \ref{Section3}, we introduce graphs of groups and their fundamental groups, as a key step in answering the realisability question in the context of pregroups. Indeed, on the one hand we have Rimlinger's results that we mentioned above (see Theorems \ref{thm rimlinger 1} and \ref{thm rimlinger 2}). On the other hand, Leary \& Stancu proved in \cite{LS} that every fusion system is realisable by a fundamental group of a graph of groups. In Section \ref{Section4}, after stating definitions and explaining how to build fusion systems from partial groups, we prove that the fusion systems induced on a Sylow $p$-subgroup by a pregroup or by its universal group are essentially the same.

\begin{thm}[Theorem \ref{thm système prégroupe} below] \label{thm:intro-F(P)=F(U(P))}
Let $P$ be a pregroup and $S$ be a finite $p$-group. Then $S$ embeds as a Sylow $p$-subgroup of $P$ if and only if it embeds as a Sylow $p$-subgroup of $U(P)$, and in this case we have $\F_S \big(U(P) \big)=\F_S(P)$. 
\end{thm}

Theorems \ref{thm:intro-universal=pi1} and \ref{thm:intro-F(P)=F(U(P))} together have the following corollary.

\begin{cor}
Let $\Lc$ be a centric linking locality associated to a fusion system $\Fc$ over a $p$-group $S$. If $\Lc$ is a pregroup, then $S$ is a Sylow $p$-subgroup of $\pi_1(B\Lc)$ and $\Fc_S(\pi_1(B\Lc))=\Fc_S(\Lc)=\Fc$.
\end{cor}

This last corollary (which can be generalised to any locality associated to a fusion system) can give some insights on the fundamental group of the classifying space of a locality associated to a fusion system. Finally, we bring all the pieces together and answer the realisability question in the context of pregroups, with the following result.

\begin{thm}[Corollary \ref{cor:realisability pregroups} below]\label{thm:intro-realisability}
Let $\Fc$ be a fusion system over a finite $p$-group $S$. There exists a finite pregroup $P$, containing $S$ as a Sylow $p$-subgroup, such that $\Fc=\Fc_S(P)$.
\end{thm}

We finish the paper by giving in Section \ref{sec:Examples} two explicit constructions of pregroups realising fusion systems. These examples are derived from the already known constructions of infinite groups realising fusion systems (Leary and Stancu \cite{LS} and Robinson \cite{Ro} respectively). For each case we also ask whether the constructed pregroup could be equipped with a structure of locality, which seems to be rarely possible.

\section{Partial groups and localities}

\label{Section1}

\subsection{Chermak's partial groups and localities}

The notions of partial groups and localities are due to Andrew Chermak. We present here the definitions and some useful properties, but more details can be found in \cite[Section 2]{Ch0} or in the preprint \cite[Section 1]{Ch1}.

For a set $X$, we denote the free monoid on $X$ by $\Wb(X)$, and for two words $u,v\in\Wb(X)$, we denote the concatenation of $u$ and $v$ by $u\circ v$. We also identify $X$ with the subset of words of length 1 in $\Wb(X)$. Finally, given two sets $X$ and $Y$ and a map $\varphi\colon X\to Y$, we will denote by $\overline{\varphi}\colon \Wb(X)\to\Wb(Y)$ the map induced by $\varphi$ defined by $\overline{\varphi}(u)=(\varphi(x_1),\varphi(x_2),\cdots,\varphi(x_n))$ for any $u=(x_1,x_2,\cdots,x_n)\in\Wb(X)$.

\begin{defi}\label{def:PG}
Let $\Mc$ be a set and let $\Db\subseteq \Wb(\Mc)$ be a subset such that,
\begin{enumerate}[label=(D\arabic*)]
\item\label{cond:D1} $\Mc\subseteq \Db$; and
\item\label{cond:D2} $u\circ v\in\Db\Rightarrow u,v\in\Db$ \quad (in particular, the empty word $\emptyset$ belongs to $\Db$).
\end{enumerate}
A mapping $\Pi:\Db\rightarrow\Mc$ is a \emph{product} if
\begin{enumerate}[label=(P\arabic*)]
\item\label{cond:P1} $\Pi$ restricts to the identity on $\Mc$; and
\item\label{cond:P2} if $u\circ v\circ w\in \Db$ then $u\circ \Pi(v)\circ w\in\Db$ and
\[\Pi(u\circ v\circ w)=\Pi\big(u\circ\Pi(v)\circ w\big).\]
\end{enumerate}
The \emph{unit} of $\Pi$ is then defined as $\Pi(\emptyset)$ and we will denote it by $1_\Mc$, or $1$ when there is no ambiguities.

An \emph{inversion} on $\Mc$ is an involutory bijection $x\mapsto x^{-1}$ on $\Mc$ together with the induced mapping $u\mapsto u^{-1}$ on $\Wb(\Mc)$ defined by,
\[u=(x_1,x_2,\cdots,x_n)\mapsto (x_n^{-1},x_{n-1}^{-1},\cdots,x_1^{-1}).\]
A \emph{partial group} is a tuple $\left(\Mc,\Db,\Pi,(-)^{-1}\right)$ where $\Pi$ is a product on $\Db$ and $(-)^{-1}$ is an inversion on $\Mc$ satisfying
\begin{enumerate}[resume*]
\item\label{cond:P3} If $u\in\Db$ then $u^{-1} \circ u \in\Db$ and $\Pi(u^{-1}\circ u)=1$.
\end{enumerate}
We will denote by $\Mc$ or $(\Mc,\Db)$ a partial group when the rest of the data is understood. The set $\Db =: \Db(\Mc)$ is called the \emph{domain} of the partial group. 
\end{defi}

Several useful properties follow easily from the above axioms.

\begin{lem} \label{lem propriétés groupes partiels}
Let $(\Mc,\Db)$ be a partial group.
\begin{enumerate}[$(1)$]
\item\label{lem partial multiplicativité} If $u \circ v \in \Db$, then $\big(\Pi(u),\Pi(v)\big) \in \Db$ and $$ \Pi(u \circ v) \ = \ \Pi(u)\Pi(v) $$ where $\Pi(u)\Pi(v)$ is short for $\Pi\big(\Pi(u),\Pi(v)\big)$.
\item\label{lem partial simplification} If $u \circ v \in \Db$, then $u^{-1} \circ u \circ v \in \Db$, $u \circ v \circ v^{-1} \in \Db$ and we have $$ \Pi(u^{-1} \circ u \circ v ) \ = \ \Pi(v) \qquad \textup{and} \qquad \Pi(u \circ v \circ v^{-1}) \ = \ \Pi(u) \ .$$
\end{enumerate}
\end{lem}

\begin{dem}
These are Lemma 2.2. (a) and (d) in \cite{Ch0}.
\end{dem}

\begin{ex}
Any group $G$ forms a partial group, setting $\Db=\Wb(G)$ and taking the product and inversion induced by group operations in $G$.
Reciprocally, if $\Mc$ is a partial group whose domain is $\Db=\Wb(\Mc)$, then $\Mc$ is a group via the binary operation $(x,y)\in\Mc^2\mapsto \Pi(x,y)\in\Mc$. 
\end{ex}

\begin{ex}
Let $\Fb(a)=\{1,a,a^{-1}\}$. We define the \emph{non-degenerated} words of $\Db_a$ to be all possible words in $\Wb(\Fb(a))$ formed by alternating $a$ and $a^{-1}$. Equivalently, the non-degenerated words of $\Db_a$ are all the different finite subwords of the ``infinite word'' $(a,a^{-1},a,a^{-1},a,a^{-1},\cdots)$. The inversion $(-)^{-1}$ is understood and, for any word $u\in\Db_a$,
\[
\Pi(u)=
\begin{cases}
\ 1 &\text{if the number of $a$'s equals the number of $a^{-1}$'s},\\
\ a &\text{if the number of $a$'s exceeds the number of $a^{-1}$'s (necessarily by 1)},\\
\ a^{-1} &\text{if the number of $a^{-1}$'s exceed the number of $a$'s (necessarily by 1)}.
\end{cases}
\]
One can then check that $\left(\Fb(a),\Db_a,\Pi_a,(-)^{-1}\right)$ defines a partial group.
\end{ex}

This last example is actually the \emph{free partial group on the set $\{a\}$} as detailed in \cite[Lemma 1.12]{Ch1}.

\medskip

Together with the notion of partial group come those of partial subgroup and morphism of partial groups.

\begin{defi}
Let $\left(\Mc,\Db,\Pi,(-)^{-1}\right)$ be a partial group. A \emph{partial subgroup} of $\Mc$ is a subset $\Nc \subseteq \Mc$ such that $\left(\Nc,\Db\cap\Wb(\Nc),\Pi_{\vert \Db\cap\Wb(\Nc)},(-)^{-1}\right)$ is a partial group. 
If $\Db\cap\Wb(\Nc)=\Wb(\Nc)$, we say that $\Nc$ is a \emph{subgroup} of $\Mc$, and if its order is a power of some prime number $p$, we say that it is a \emph{$p$-subgroup} of $\Mc$.
\end{defi}

\begin{defi}
Let $\left(\Mc_1,\Db_1,\Pi_1,(-)^{-1}\right)$ and $\left(\Mc_2,\Db_2,\Pi_2,(-)^{-1}\right)$ be two partial groups.
A \emph{morphism of partial groups}, or \emph{partial group homomorphism}, is a map $\varphi \colon \Mc_1\to \Mc_2$ such that 
\begin{enumerate}[label=(H\arabic*)]
\item $\overline{\varphi}\left(\Db_{1}\right)\subseteq \Db_{2}$;
\item for any $u\in\Db_1$, $\Pi_2\left(\overline{\varphi}(u)\right)=\varphi\left(\Pi_1(u)\right)$.
\end{enumerate}
Moreover, $\varphi \colon \Mc_1\to \Mc_2$ is called an \emph{isomorphism} of partial groups if the map $\varphi$ is bijective and if $\varphi^{-1}$ is also a morphism of partial group. Finally, an \emph{automorphism} of $\Mc_1$ is an isomorphism $\varphi\colon \Mc_1\to\Mc_1$.
\end{defi}

With this notion of morphisms and the usual composition on maps, the class of partial groups forms a category $\ParG$ which contains the category of groups as a full subcategory.

Chermak introduced partial groups to study the $p$-local structure of finite groups when $p$ is a prime number. For that purpose he defined the notion of \emph{locality}, which allows to encode and manipulate these $p$-local structures. In order to define localities, we first need to talk about conjugation in partial groups, and define \emph{objective partial groups}.

\begin{nota}\label{notation conjugaison}
Given $\Mc$ a partial group and $g \in \Mc$, we denote by $\Db(g)$ the set of all $x \in \Mc$ such that $(g^{-1},x,g) \in \Db$, and by $c_g \colon x \mapsto x^g$ the map sending $x \in \Db(g)$ to $\Pi(g^{-1},x,g)$. If in addition $X$ is a subgroup of $\Mc$ such that $X \subseteq \Db(g)$, then we denote by $X^g$ the set of all $x^g$ for $x \in X$. Beware that $X^g$ is not a subgroup of $\Mc$ in general.
\end{nota}

\begin{lem} \label{lem partial conjugaison}
Let $\Mc$ be a partial group and $g \in \Mc$. Then $c_g$ defines a bijection $\Db(g) \to \Db(g^{-1})$ whose inverse is $c_{g^{-1}}$.
\end{lem}

\begin{dem}
This is Lemma 2.5. (c) in \cite{Ch0}.
\end{dem}

\begin{defi} \label{def objective}
Let $\Mc$ be a partial group and $\Delta$ a collection of subgroups of $\Mc$.
Define $\Db_\Delta$ to be the set of all $w=(g_1,g_2,\cdots,g_n)\in\Wb(\Mc)$ such that
\begin{equation}\tag{\ding{92}}\label{*}
\exists (X_0,X_1,\cdots,X_n)\in\Wb(\Delta), \ \forall i \in \{1, \cdots ,n\}, \qquad X_{i-1} \subseteq \Db(g_i) \quad \textup{and} \quad (X_{i-1})^{g_i}=X_i \ .
\end{equation}
We then say that $(\Mc,\Delta)$ is an \emph{objective partial group} if the following two conditions holds.
\begin{enumerate}[label=(O\alph*)]
\item $\Db(\Mc)=\Db_\Delta$.\label{O1}
\item Whenever $X$ and $Y$ are in $\Delta$ and $g\in\Mc$ is such that $X \subseteq \Db(g)$ and $X^g$ is a subgroup of $Y$, then every subgroup of $Y$ containing $X^g$ is in $\Delta$.\label{O2}
\end{enumerate}
\end{defi}

Hence an objective partial group is a particular instance of partial group whose domain is given by the composable conjugation maps between a fixed set of subgroups. For example, given any group $G$ and any collection $\Delta$ of subgroups of $G$, then $(G,\Db_\Delta))$ is an objective partial group.

\begin{defi} \label{def locality}
Let $p$ be a prime number, let $\Lc$ be a finite partial group. Let $S$ be a $p$-subgroup of $\Lc$, and let $\Delta$ be a collection of subgroups of $S$ such that $S\in\Delta$. We say that $(\Lc,\Delta,S)$ is a \emph{locality} if:
\begin{enumerate}[label=(L\alph*)]
\item $(\Lc,\Delta)$ is objective; and
\item $S$ is maximal in the poset (ordered by inclusion) of finite $p$-subgroups of $\Lc$.
\end{enumerate}
\end{defi}

%
%

\subsection{A simplicial point of view on partial groups}

Let $\left(\Mc,\Db,\Pi,(-)^{-1}\right)$ be a partial group.  Broto and Gonzales \cite{BG} pointed out that $\Pi$ induces a simplicial set structure on $\Db$ in the same way that the product in a group $G$ induces a simplicial set structure on $\Wb(G)$. 

\begin{defi}
Let $\left(\Mc,\Db,\Pi,(-)^{-1}\right)$ be a partial group. We denote by $\Bc(\Mc)$ the simplicial set whose $n$-simplices, for $n \in \Nb$, are the elements of $\Bc_n(\Mc)=\Db_n$, the set of words of length $n$ in $\Db$. The face operators are given, for $n \in \Nb^*$, $i\in\{0,1,\cdots,n-1\}$ and $(m_1,m_2,\cdots,m_n)\in \Bc_n(\Mc)$, by
 
 \[d_i(m_1,m_2,\cdots,m_n)=
  \begin{cases}
     \  (m_2,m_3,\cdots,m_n) &\text{if }i=0 ,\\
     \  (m_1,\cdots,\Pi(m_{i},m_{i+1}),\cdots,m_n) & \text{if }1\leq i\leq n-1 ,\\
     \  (m_1,m_2,\cdots,m_{n-1})&\text{if }i=n ;
  \end{cases}
 \]
 and the degeneracy operators are defined, for $n \in \Nb$, $i\in\{0,1,\cdots,n\}$ and $(m_1,m_2,\cdots,m_n)\in \Bc_n(\Mc)$, by
 \[s_i(m_1,m_2,\cdots,m_n)=(m_1,\cdots,m_i,1,m_{i+1},\cdots,m_n)\]
 where 1 is the unit of $\Mc$.
 
 Finally, the geometric realisation of this simplicial set will be denoted by $B\Mc:=|\Bc(\Mc)|$ and called the \emph{classifying space} of $\Mc$.
\end{defi}

As mentioned above, when we are working with an actual group $G$ (i.e. $\Db(G)=\Wb(G)$) then $\Bc(G)$ is the classical bar construction and $BG$ is a classifying space for the group $G$. 

Notice also that a map between two partial groups $f\colon \Mc_1\to \Mc_2$ is a partial group homomorphism if and only if the map induced by $f$ on words, $\overline{f}\colon \Bc(\Mc_1)\to\Bc(\Mc_2)$, is a simplicial map.  Broto and Gonzalez actually showed that this constitutes a fully faithfull embedding $\Bc\colon \ParG\to \sSet$ of $\ParG$ into $\sSet$ the category of simplicial sets (see \cite[Section 2]{BG} for more details).

There is a deep connection between a partial group and the geometric realisation of its associated simplicial set. For example, there is a correspondence between extension of partial groups and fiber bundles of the corresponding simplicial sets as highlighted by Broto and Gonzalez \cite{BG,Go}. Here we focus on the fundamental group $\pi_1(B\Mc)$ of the geometric realisation. 

\begin{prop}\label{prop:fundamental group}
Let $\Mc$ be a partial group. There is a natural morphism of partial groups $\theta\colon \Mc \to \pi_1(B\Mc)$ and for any 
group $G$ and any morphism of partial groups $\alpha\colon \Mc\to G$, there is a unique group homomorphism $\overline{\alpha}\colon \pi_1(B\Mc)\to G$ such that $\alpha= \overline{\alpha} \circ \theta$. 
\end{prop}

\begin{dem}
This is a direct consequence of \cite[Prop III.2.7]{AKO}. If we denote $K:=\Bc(\Mc)=\Db(\Mc)$ to follow the notation of \cite[Prop III.2.7]{AKO}, we have that $K_0$ is just a point given by the empty word (and it will be denoted by $x_0$ in the rest of the proof), and $K_1=\Mc$. Then the map $\theta\colon K_1\to\pi_1(|K|,x_0)$ of \cite[Prop III.2.7]{AKO} is just a map $\theta\colon \Mc\to \pi_1(|K|,x_0)$ such that for all $(g,h)\in \Db_2(\Mc)=K_2$, we have $\theta\left(\Pi(gh)\right)=\theta(g)\theta(h)$. This is equivalent to state that $\theta$ is a partial group homomorphism from $\Mc$ to the group $\pi_1(B\Mc)$. Finally, the universal property is just a restatement of the second part of \cite[Prop III.2.7]{AKO}.
\end{dem}

In particular, $B\colon \ParG\to\Top_*$ is a functor from $\ParG$ to the category $\Top_*$ of pointed topological spaces (because if $\Mc$ is a partial group, $\Bc(\Mc)$ has only one 0-simplex, so there is a canonical choice of basepoint for $B\Mc$), and if we denote the category of groups by $\Grps$, the functor $\ParG\to \Grps$ defined as the composition of $B$ followed by the fundamental group functor $\pi_1\colon \Top_*\to\Grps$, is left adjoint to the forgetful functor $U\colon \Grps\to \ParG$.

\section{Stallings' pregroups and their universal groups}

\label{Section2}

\subsection{Generalities on pregroups}

The notion of pregroup was introduced by Stallings in \cite{St}. This generalisation of the group structure aims at providing ``nice'' generating sets for certain families of groups, such as amalgamated products of groups. Most of the material here comes from \cite[Section 3.A]{St} and \cite[Part I, section I]{Ri}.

\begin{defi}\label{def:Pre}
A \emph{pregroup} is a tuple $\big(P,D,m,1_P,(-)^{-1} \big)$ where $P$ is a set, $D$ is a subset of $P\times P$ called the \emph{domain}, $m\colon D\to P$ and $(-)^{-1} \colon P\to P$ are maps and $1_P$ is an element of $P$ called the \emph{unit}, such that for all $w,x,y,z\in P$ we have:
\begin{enumerate}[label = (Pr\arabic*)]
 \item\label{cond:Pre1} $(1_P,x)\in D$, $(x,1_P)\in D$ and $m(1_P,x)=m(x,1_P)=x$.
 \item\label{cond:Pre2} $(x,x^{-1})\in D$, $(x^{-1},x)\in D$ and $m(x,x^{-1})=m(x^{-1},x)=1_P$.
 \item\label{cond:Pre3} If $(w,x),(x,y) \in D$, then 
 \[ \big(m(w,x),y \big) \in D \quad \Longleftrightarrow \quad \big(w,m(x,y)\big) \in D \]
 and in that case, $m \big(m(w,x),y \big) = m\big(w,m(x,y)\big)$.
 \item\label{cond:Pre4} If $(w,x),(x,y),(y,z) \in D$, then $\big(w,m(x,y) \big)\in D$ or $\big(m(x,y),z \big)\in D$.
\end{enumerate}
To simplify the reading, when $(x,y) \in D$, we will denote $xy$ instead of $m(x,y)$. When $w,x,y \in P$ satisfy $(w,x),(x,y) \in D$ and $(wx,y)\in D$ (or equivalently, $(w,xy)\in D$ by \ref{cond:Pre3}), we say that $(w,x,y)$ \emph{associates}.
\end{defi}

Considering the above definition, a pregroup $\big(P,D,m,1_P,(-)^{-1} \big)$ is a group if and only if $D = P \times P$. Moreover, we deduce easily the following properties from the axioms.
\begin{lem} \label{lem propriétés prégroupe}
Let $\big(P,D,m,1_P,(-)^{-1} \big)$ be a pregroup. We have the following:
\begin{enumerate}[$(1)$]
\item\label{lem item involution} If $x \in P$, then $\big(x^{-1}\big)^{-1} = x$.
\item\label{lem item simplification} If $(x,y) \in D$, then $(x^{-1},xy) \in D$ and $x^{-1}(xy) = y$. Similarly, $(xy,y^{-1}) \in D$ and $(xy)y^{-1} = x$.
\item\label{lem item inverse} Let $x,y \in P$. Then $(x,y) \in D$ if and only if $(y^{-1},x^{-1}) \in D$, and in this case we have $y^{-1}x^{-1} = (xy)^{-1}$.
\item\label{lem item intercalation} Let $a \in P$ such that $(x,a),(a^{-1},y) \in D$. Then $(x,y) \in D$ if and only if $(xa,a^{-1}y) \in D$. In this case, we have $(xa)(a^{-1}y) = xy$.
\end{enumerate}
\end{lem}

\begin{dem}
These properties are proved in \cite{LiShi} for \ref{lem item inverse}, and in \cite{St} for the others. However we prefer to give them a proof here to ensure consistency (see next remark).

Let us start with \ref{lem item involution}. By Axiom \ref{cond:Pre2}, $(x,x^{-1})$ and $\big(x^{-1},(x^{-1})^{-1}\big)$ are in $D$, and $xx^{-1} = 1_P = x^{-1}(x^{-1})^{-1}$. By Axiom \ref{cond:Pre1}, $\big(xx^{-1},(x^{-1})^{-1}\big)$ and $\big(x,x^{-1}(x^{-1})^{-1}\big)$ are in $D$ and their products equal $(x^{-1})^{-1}$ and $x$ respectively. Finally, by Axiom \ref{cond:Pre3}, these two products are equal.

For \ref{lem item simplification}, by Axiom \ref{cond:Pre2} we have $(x^{-1},x) \in D$ and $x^{-1}x = 1_P$. Thus by Axiom \ref{cond:Pre1}, $(x^{-1}x , y) \in D$ and $(x^{-1}x)y = y$, so $(x^{-1},xy) \in D$ and $x^{-1}(xy) = (x^{-1}x)y = y$ by Axiom \ref{cond:Pre3}. The other case of \ref{lem item simplification} follows from the same arguments.

Now for \ref{lem item inverse}, by \ref{lem item involution} it is enough to prove that $(x,y) \in D$ implies $(y^{-1},x^{-1}) \in D$. We know that $\big( (xy)^{-1} , (xy)y^{-1} \big) \in D$ because of \ref{lem item simplification}, and $\big((xy)y^{-1},x^{-1} \big) = (x,x^{-1}) \in D$ by Axiom \ref{cond:Pre2}. Moreover $\Big( (xy)^{-1}, \big((xy)y^{-1}\big)x^{-1} \Big) = \big( (xy)^{-1} , 1_P \big) \in D$. Hence using Axiom \ref{cond:Pre3} we deduce that $\Big( (xy)^{-1}\big((xy)y^{-1}\big),x^{-1} \Big) = \big( y^{-1}, x^{-1} \big) \in D$ and $y^{-1}x^{-1} = (xy)^{-1}$.

Finally for \ref{lem item intercalation}, we know from \ref{lem item simplification} that $(xa,a^{-1}) \in D$. Since we also have $(a^{-1},y) \in D$, the result becomes nothing more than Axiom \ref{cond:Pre3}.
\end{dem}

\begin{rem}
In \cite{St}, Stallings introduces the definition of pregroup with property \ref{lem item inverse} from Lemma \ref{lem propriétés prégroupe} as an extra axiom. It appeared later to be redundant, so that this property is no more included as an axiom in \cite{Ri}. On the contrary, in the latter Rimlinger assumes that $x \mapsto x^{-1}$ is an involution by definition, but this is a consequence of the other axioms, as proved in \cite{St}.
\end{rem}

As for partial groups, we have the notions of subpregroup and morphism of pregroups.

\begin{defi}
Let $\big(P,D,m,1_P,(-)^{-1} \big)$ be a pregroup. A subset $Q\subset P$ induces a \emph{subpregroup} of $P$ if the tuple $\big(Q, D\cap (Q\times Q), m|_{D\cap (Q\times Q)},1_P,(-)^{-1} \big)$ is a pregroup. If $Q\times Q\subseteq D$ we will call it a \emph{subgroup} of $P$. 
\end{defi}

\begin{defi}
Let $\big(P_1,D_1,m_1,1_1,(-)^{-1} \big)$ and $\big(P_2,D_2,m_2,1_2,(-)^{-1} \big)$ be two pregroups. A \emph{morphism of pregroups}, or \emph{pregroup homomorphism}, is a map $\varphi\colon P_1\to P_2$ such that for every $(x,y)\in D_1$, $(\varphi(x),\varphi(y))\in D_2$ and $m_2(\varphi(x),\varphi(y))=\varphi(m_1(x,y))$.
\end{defi}

Thus pregroups, together with morphisms of pregroups, define a category that we will denote by $\PrG$. One important fact about pregroups is that we can associate to any of them a group called its \emph{universal group}.

\begin{defi}
Let $\big(P,D,m,1_P,(-)^{-1} \big)$ be a pregroup. The \emph{universal group} of $P$, denoted by $U(P)$, is the group with presentation 
\[
U(P)=\langle \ P \ \vert \ m(x,y)y^{-1}x^{-1} \textup{ for all } (x,y)\in D \ \rangle
\]
where $m(x,y)y^{-1}x^{-1}$ is the product of $m(x,y)$, $y^{-1}$ and $x^{-1}$ in the free group generated by $P$ (not in $P$).
\end{defi}

\begin{nota}
In the rest of the paper, a pregroup $\big(P,D,m,1_P,(-)^{-1} \big)$ will be often just denoted by $(P,D)$, or even $P$ when reference to $D$ is not needed.
\end{nota}

Before continuing with properties of the universal group, let us introduce a family of examples of pregroups, naturally arising from amalgamated products of groups.

\begin{ex}[{\cite[Example 3.A.5.1]{St}}]\label{ex:A*_CB}
Let $A,B$ and $C$ be three groups, and let $\varphi_A\colon C\to A$ and $\varphi_{B}\colon C\to B$ be injective group homomorphisms. 
Set $A\cup_C B := A\sqcup B/\sim$, where $\sim$ is defined by $\varphi_A(x) \sim \varphi_B(x)$ for all $x \in C$. We can identify $A \cup_C B$ with $(A\setminus \varphi_A(C)) \sqcup C \sqcup (B \setminus \varphi_B(C))$, $A$ with $(A\setminus \varphi_A(C)) \sqcup C$ and $B$ with $(B\setminus \varphi_B(C)) \sqcup C$, so that $A \cup_C B$ contains both $A$ and $B$, and $A \cap B = C$ with these identifications. Then $A \cup_C B$ has a natural pregroup structure with domain $$ D := \left\{ \ (x,y) \in (A \cup_C B) \times (A \cup_C B) \ \vert \ x,y \in A \ \textup{ or } \ x,y \in B \ \right\} $$ and the obvious operations. In this way it can be seen as a subset of the amalgamated product $A *_C B$ (see example (1) in \ref{examples graphs of groups} for the definition), and in fact one can easily see that $U(A\cup_C B)=A *_C B$.
\end{ex}

The following universal property is a direct consequence of the definition of $U(P)$.

\begin{prop}\label{prop:universalgroup}
Let $P$ be a pregroup. 
The natural map $\iota_P \colon P \to U(P)$ is a morphism of pregroups, and for any group $G$ together with a pregroup homomorphism $\alpha\colon P\to G$, there exists a unique group homomorphism $\overline{\alpha}\colon U(P)\to G$ such that $\alpha=\overline{\alpha}\circ \iota_P$.
\[
\xymatrix{ U(P) \ar@{-->}^{\exists!\;\overline{\alpha}}[rr]& & G\\ P\ar^{\iota_P}[u]\ar_{\alpha}[rru]&&}
\] 
\end{prop}

In Theorem 3.A.4.5 from \cite{St}, Stallings proved that the universal group of a pregroup has a solvable word problem. A noteworthy corollary of this theorem is that $P$ is injectively embedded in $U(P)$.

\begin{thm}[{\cite[Corollary 3.A.4.6]{St}}]
The morphism $\iota_P \colon P\to U(P)$ is an injective pregroup homomorphism.
\end{thm}

\begin{nota}
Let $P$ be a pregroup. In general we will identify $P$ with its image under $\iota_P \colon P\to U(P)$, and for $(x_1,x_2,\cdots,x_n)\in \Wb(P)$ we will denote by $x_1x_2\cdots x_n$ its product in $U(P)$. 
\end{nota}

We will also need a weaker version of Stallings' theorem, which we state below, just after defining $P$-reduced words.

\begin{defi}
A word $(x_1,x_2,\cdots,x_n)\in \Wb(P)$ is said to be \emph{$P$-reduced} if for all $i\in \{1,2,\cdots,n-1\}$, $(x_i,x_{i+1})\not\in D$. The empty word is $P$-reduced.
\end{defi}

\begin{thm}[{\cite[Theorem 3.A.4.5]{St}}]\label{th:wordsolv}
Let $P$ be a pregroup and let $(x_1,x_2,\cdots,x_n)\in \Wb(P)$ and $(y_1,y_2,\cdots,y_m)\in\Wb(P)$ be two $P$-reduced words.
If $x_1x_2\cdots x_n=y_1y_2\cdots y_m$ then $n=m$. 
\end{thm}

In particular, for $n=2$, this result provides a characterisation of the domain $D$ in terms of products in the universal group.

\begin{cor}\label{cor:domain=product}
Let $(P,D)$ be a pregroup and $x,y\in P$.
Then $(x,y)\in D$ if and only if $xy \in P$, where the product $xy$ is performed in $U(P)$.
\end{cor}

We also get the following corollary.

\begin{cor} \label{cor:subgroups of pregroup}
Let $P$ be a pregroup and $U(P)$ be its universal group. Given any subgroup $H$ of $U(P)$, $H$ is a subgroup of $P$ if and only it is included in $P$.
\end{cor}

\subsection{Elements of finite order, finite subgroups and conjugation in the universal group}

In this section, $(P,D)$ will be a fixed pregroup and $U(P)$ its universal group. We will have a look at elements of finite order in the universal group.

\begin{defi}
A \emph{cyclic element} of $P$ is an element $x\in P$ such that $\langle x\rangle$ is a subgroup of $P$ (and not just of $U(P)$). In view of Corollary \ref{cor:subgroups of pregroup}, this is equivalent to ask for $\langle x \rangle$ to be included in $P$.
\end{defi}

\begin{lem}\label{lem:cyclic elements}
Let $x$ be an element of $P$. Then the following are equivalent:
\begin{enumerate}[$(1)$]
\item \label{lemcyclic1} $x$ is cyclic,
\item \label{lemcyclic2} $(x,x) \in D$,
\item \label{lemcyclic3} $x^2 \in P$.
\end{enumerate}
\end{lem}

\begin{dem}
Rimlinger already proved this result in \cite[Corollary 1.10]{Ri}, but the proof is short and is an easy (warm-up) example for using Axiom \ref{cond:Pre4} in a proof involving pregroups, so we give it here.

The equivalence between \ref{lemcyclic2} and \ref{lemcyclic3} is a direct application of Corollary \ref{cor:domain=product}, and the fact that \ref{lemcyclic1} implies \ref{lemcyclic3} is trivial. Now we assume \ref{lemcyclic2} and we prove by induction on $n \in \Nb$ that $x^n \in P$. The cases $n=0$ and $n=1$ are trivial. Let $n \geq 1$ and suppose that $x^i \in P$ for all $i \in \{1, \cdots ,n \}$. In particular we have $(x^{n-1},x) \in D$, $(x,x^{n-1}) \in D$ by induction hypothesis and Corollary \ref{cor:domain=product}, and $(x,x) \in D$ by \ref{lemcyclic2}. Applying Axiom \ref{cond:Pre4} to the tuple $(x,x^{n-1},x,x)$, we get $(x,x^n) \in D$ or $(x^n,x) \in D$, and both cases give $x^{n+1} \in P$ by a final application of Corollary \ref{cor:domain=product}. This concludes the proof by induction.
\end{dem}

\begin{lem}\label{lem:cyclic conjugate}
Let $x\in U(P)$ be an element of finite order. 
\begin{enumerate}[$(i)$]
\item \label{lemi} If $x\in P$ then $x$ is a cyclic element of $P$.
\item \label{lemii} If $x\not\in P$, then $x$ is conjugate (in $U(P)$) to a cyclic element of $P$.
\end{enumerate}
\end{lem}

\begin{dem}
By hypothesis, there exists $r \in \Nb^*$ such that $x^r=1$ in $U(P)$. If $x \in P$, then by Theorem \ref{th:wordsolv}, the word $(x,x,\cdots,x)\in P^r$ cannot be $P$-reduced. Therefore $(x,x)\in D$ and, by Lemma \ref{lem:cyclic elements}, $x$ is a cyclic element of $P$. This proves \ref{lemi}.

\medskip

Thanks to \ref{lemi}, and since the conjugate of an element of finite order is of finite order, it is enough to show that $x$ is conjugate to an element of $P$ to prove \ref{lemii}.

Let $w=(p_1,p_2,\cdots,p_n)\in \Wb(P)$ be a minimal $P$-reduced word such that $y=p_1p_2\cdots p_n$ is conjugate to $x$ in $U(P)$ and assume that $n\geq 2$. Then $y$ is also an element of finite order so there exists $r \in \Nb^*$ such that $y^r=1$. This implies that the concatenation $w^r$ of $w$ with itself $r$ times is not $P$-reduced. Since $w$ is $P$-reduced, this implies that $(p_n,p_1)\in D$. Therefore $p_1^{-1}yp_1=p_2p_3\cdots p_{n}p_1$ is conjugate to $x$ and is the product of the word $(p_2,p_3,\cdots,p_np_1)$ which is of length $n-1$. This contradicts the minimality of $w$. Hence $n=1$, i.e. $y \in P$.
\end{dem}

In particular, this lemma allows us to talk about elements of $P$ of finite order without ambiguity.

\begin{lem}\label{lem:Prconj}
Let $x,g\in P$. If $x$ is an element of $P$ of finite order, then the following are equivalent.
 \begin{enumerate}[$(1)$]
 \item $g^{-1}x \in P$,
 \item $xg \in P$.
 \end{enumerate}
\end{lem}

\begin{dem}
According to Lemma \ref{lem:cyclic conjugate}, $x$ is cyclic. Assume that $g^{-1}x \in P$ but $xg \not\in P$. We have $x^{-1}g \in P$ by items \ref{lem item involution} and \ref{lem item inverse} of Lemma \ref{lem propriétés prégroupe}, and $(g^{-1},x), (x,x), (x,x^{-1}g) \in D$ by Lemma \ref{lem propriétés prégroupe}. Hence, by \ref{cond:Pre4}, since $xg \not\in P$, $g^{-1}x^2\in P$. Now we proceed by induction to prove that $g^{-1}x^k \in P$ for every $k \in \Nb$. We already proved it for $k \leq 2$, so let $k \in \Nb$, $k \geq 2$ and assume that $g^{-1}x^k \in P$. Then we also have $x^{-k}g \in P$ by Lemma \ref{lem propriétés prégroupe}, so that $(g^{-1},x), (x,x^k), (x^k,x^{-k}g) \in D$. Applying Axiom \ref{cond:Pre4}, we get $g^{-1}x^{k+1} \in P$, and this concludes the induction. But $x$ is assumed to be of finite order, so $x^{-1} = x^k$ for some $k \in \Nb$. Therefore, $g^{-1}x^{-1} \in P$, so $xg \in P$ by item \ref{lem item inverse} of Lemma \ref{lem propriétés prégroupe}, which contradicts the initial assumption. 

For the other implication, notice that $x^{-1}$ is also an element of $P$ of finite order, and $xg \in P$ if and only if $g^{-1} x^{-1}\in P$. Thus, applying the previous implication, we get $x^{-1}g \in P$, which is equivalent to $g^{-1}x \in P$.
\end{dem}

Now we establish some facts about conjugation in $U(P)$, beginning with a technical lemma.

\begin{lem} \label{lem:technical conjugation}
Let $k \in \Nb^*$, let $x$ be an element of $P$ of finite order, and let $(g_0,g_1,\cdots, g_{k-1})\in\Wb(P)$ be a $P$-reduced word such that $g_{k-1}^{-1}g_{k-2}^{-1}\cdots g_0^{-1}xg_0g_{1}\cdots g_{k-1} \in P$. Then for every $i \in \{1,2,\cdots, k\}$, $g_{i-2}^{-1}\cdots g_0^{-1}xg_0g_{1} \cdots g_{i-2} \in P$ and $\left(g_{i-1}^{-1},g_{i-2}^{-1}\cdots g_0^{-1}xg_0g_{1}\cdots g_{i-2},g_{i-1} \right)$ associates.
\end{lem}

\begin{dem}
We proceed by induction on $k$. For $k=1$, let $g\in P$ and $x$ be a cyclic element of $P$ such that $gxg^{-1}\in P$. In particular, we know that the word $(g^{-1},x,g)$ is not $P$-reduced thanks to Theorem \ref{th:wordsolv}. Thus $g^{-1}x \in P$ or $xg \in P$, and by Lemma \ref{lem:Prconj} this implies that $g^{-1}x \in P$ and $xg \in P$. Since $(g^{-1}x)g$ is in $P$, $(g^{-1},x,g)$ associates. 

Now let $k \in \Nb^*$ and assume the result is true for the rank $k$. Let $(g_0,g_1,\cdots, g_k)\in\Wb(P)$ be a $P$-reduced word, let $x$ be a cyclic element of $P$ and assume that $g_{k}^{-1}g_{k-1}^{-1}\cdots g_0^{-1}xg_0g_{1}\cdots g_{k} \in P$. In particular, the word $(g_{k}^{-1},g_{k-1}^{-1},\cdots, g_0^{-1},x,g_0,g_{1},\cdots, g_{k})$ is not $P$-reduced according to Theorem \ref{th:wordsolv}. Since $(g_0,g_1,\cdots, g_k)$ is $P$-reduced, this is also the case of $(g_{k}^{-1},g_{k-1}^{-1},\cdots, g_0^{-1})$, so we have $g_0^{-1}x\in P$ or $xg_0 \in P$. By Lemma \ref{lem:Prconj}, they are both in $P$. Moreover, the word $(g_{k}^{-1},g_{k-1}^{-1},\cdots,g_1^{-1}, g_0^{-1}x,g_0,g_{1},\cdots, g_{k})$ is still not $P$-reduced. Therefore $g_1^{-1}(g_0^{-1}x)\in P$ or $(g_0^{-1}x)g_0 \in P$. Assume the latter is true, then we can conclude directly using the induction hypothesis, because the conjugate of an element of finite order is again of finite order. Else, if $g_1^{-1}(g_0^{-1}x)\in P$ then $(g_1^{-1},g_0^{-1}x),(g_0^{-1}x,x^{-1}),(x^{-1},xg_0) \in D$ and by \ref{cond:Pre4}, this implies $g_0^{-1}xg_0=(g_0^{-1}xx^{-1})xg_0\in P$ or $g_1^{-1}g_0^{-1}=g_1^{-1}(g_0^{-1}xx^{-1})\in P$ (which is absurd). Thus we also get $g_0^{-1}xg_0 \in P$, and we can apply the induction hypothesis, taking $g_0^{-1}xg_0 \in P$ for the element of finite order and $(g_1,g_2,\cdots,g_k)$ for the $P$-reduced word.
\end{dem}

\begin{prop}\label{prop:Prconj}
Let $x,y\in P$ be two elements of finite order. If $x$ and $y$ are conjugate in $U(P)$ then there exist two finite sequences $x=:x_0,x_1,\cdots,x_{k-1},x_k:=y$ and $g_0,g_1,\cdots, g_{k-1}$ of elements of $P$ such that for all $i \in \{1, \cdots , k-1\}$, $(g_i^{-1},x_i,g_i)$ associates and $g_i^{-1}x_ig_i=x_{i+1}$.
\end{prop}

\begin{dem}
Let $g \in U(P)$ be such that $g^{-1}xg=y$, and let $(g_0,g_1,\cdots, g_{k-1})\in\Wb(P)$ be a $P$-reduced word representing $g$. Then the result follows directly from Lemma \ref{lem:technical conjugation}, if we define the $x_i \in P$ inductively by $x_{i+1} := g_i^{-1}x_ig_i$.
\end{dem}

\begin{prop} \label{prop conjugate sgs in P}
Let $Q$ and $R$ be two finite subgroups of $P$. If $Q$ and $R$ are conjugate in $U(P)$ then there exist a sequence $Q_0=Q,Q_1,\cdots,Q_k=R$ of finite subgroups of $P$ and a sequence $g_0,g_1,\cdots, g_{k-1}$ of elements of $P$, such that for each $i \in \{1, \cdots , k-1\}$ and every $x \in Q_i$, $(g_i^{-1},x,g_i)$ associates, and $g_i^{-1} Q_i g_i=Q_{i+1}$.
\end{prop}

\begin{dem}
Let $g \in U(P)$ be such that $g^{-1} Q g=R$, and let $(g_0,g_1,\cdots, g_{k-1})\in\Wb(P)$ be a $P$-reduced word representing $g$. Define the subgroups $Q_i$ of $U(P)$ inductively by $Q_0:=Q$ and $Q_{i+1} := g_i^{-1} Q_i g_i$, so that $Q_k=R$. Fix $i \in \{0, \cdots , k-1 \}$ and $x_i \in Q_i$. By definition of $Q_i$, there exists $x_0 \in Q$ such that $x_i = g_{i-1}^{-1} \cdots g_0^{-1} x_0 g_0 \cdots g_{i-1}$, where the products are performed in $U(P)$. Moreover, $g_{k-1}^{-1} \cdots g_0^{-1} x_0 g_0 \cdots g_{k-1} \in R \subseteq P$ by hypothesis. Applying Lemma \ref{lem:technical conjugation}, we get in particular that $x_i \in P$ and $(g_i^{-1},x_i,g_i)$ associates. As this holds for general $x_i \in Q_i$, we deduce that $Q_i \subseteq P$, so by Corollary \ref{cor:subgroups of pregroup}, $Q_i$ is a subgroup of $P$ and the proof is complete.
\end{dem}

\subsection{Pregroups are partial groups}

Let $\left(P,D\right)$ be a pregroup and let $U(P)$ be its universal group. Consider the following subset of $\Wb(P)$: $$\Db_P \ = \ \left\{ \ (x_1,x_2,\dots,x_n) \in \Wb(P) \ \vert \ \forall i,j\in\{1,2,\dots, n\} \textup{ with } i<j, \ \ x_ix_{i+1}\cdots x_j\in P \ \right\} $$

By repeated applications of Corollary \ref{cor:domain=product} and Axiom \ref{cond:Pre3}, there exists a well-defined map $\Pi_{P}$ from $\Db_P$ to $P$ assigning $w=(x_1,x_2,\dots,x_n)\in\Db_P$ to
\[\Pi_P(w)=x_1x_2\cdots x_n.\]

This defines a partial group structure on $P$, as stated in the following proposition.
\begin{prop}\label{prop pregroup is partial group}
Let $P$ be a pregroup. Then $\left(P,\Db_P,\Pi_P,(-)^{-1}\right)$ is a partial group. This construction induces a fully faithful functor $\PrG\to\ParG$.
\end{prop}

\begin{dem}
First, let us check the axioms for partial groups. By construction, the domain $\Db_P$ satisfies \ref{cond:D1} and \ref{cond:D2}. It is also clear by definition that $\Pi_P$ restricts to identity on $P$. Condition \ref{cond:P2} is a consequence of the fact that $\Pi_P$ is well-defined from the above formula. Finally, the fact that $x \mapsto x^{-1}$ is an involutory bijection comes from item \ref{lem item involution} in Lemma \ref{lem propriétés prégroupe}, and \ref{cond:P3} comes from Axiom \ref{cond:Pre2} for pregroups.

Now, if $P$ and $Q$ are two pregroups and $\varphi \colon P \to Q$ is a pregroup homomorphism, one has to check that $\varphi$ is in fact a partial group homomorphism between the associated partial groups. If $(x_1,x_2,\dots,x_n) \in \Db_P$, we prove by induction on $j-i$ that $\varphi(x_i) \cdots \varphi(x_j) \in Q$ and $\Pi_Q\big( \varphi(x_i), \dots, \varphi(x_j) \big) = \varphi \big( \Pi_P(x_i,\dots,x_j) \big)$ for all $i,j\in\{1,2,\dots, n\}$ such that $i<j$. For $j-i=1$, this is just the fact that $\varphi$ is a pregroup homomorphism. For $j-i>1$, by induction hypothesis we have $$ \varphi(x_i) \cdots \varphi(x_j) \ = \ \Pi_Q\big( \varphi(x_i), \dots, \varphi(x_{j-1}) \big) \cdot \varphi(x_j) \ = \ \varphi \big( \Pi_P(x_i,\dots,x_{j-1}) \big) \cdot \varphi(x_j) \ .$$ Since $\Pi_P(x_i,\dots,x_{j-1}) \cdot x_j = x_i \cdots x_j \in P$, applying the definition of pregroup homomorphism for $\varphi$, we deduce from the above that $\varphi(x_i) \cdots \varphi(x_j) \in Q$ and \begin{align*}
 \Pi_Q\big( \varphi(x_i), \dots, \varphi(x_j) \big) \ = \ \varphi(x_i) \cdots \varphi(x_j) \ &= \ \varphi \big( \Pi_P(x_i,\dots,x_{j-1}) \big) \cdot \varphi(x_j) \\
& = \ \varphi \big( \Pi_P(x_i,\dots,x_{j-1}) \cdot x_j \big) \\
& = \ \varphi \big( \Pi_P(x_i,\dots,x_{j}) \big) .
\end{align*}

Hence we get a functor $\PrG\to\ParG$, and there only remains to check that it is full and faithful. The faithfulness is obvious, since the underlying set map of a pregroup homomorphism is the same as that of the associated partial group homomorphism. Finally, the axioms for pregroup homomorphisms are just restrictions of those for partial group homomorphisms (applying them only for couples of elements in the domain), so that any partial group homomorphism between pregroups is clearly a pregroup homomorphism. This gives the fullness.
\end{dem}

\begin{ex}\label{ex:A*_CB revisted}
Let $A,B,C$ be as in Example \ref{ex:A*_CB} and $P=A\cup_C B$ be the pregroup constructed in that example. Considering the underlying set of $P$ as a subset of $A*_CB$, we then have $\Db_P=\Wb(A)\cup\Wb(B)$.
\end{ex}

Going back to the simplicial point of view on partial groups, this reveals a relation between the universal group of a pregroup and the fundamental group of the classifying space of the associated partial group.

\begin{cor}\label{cor:universal=pi1}
Let $P$ be a pregroup and let $\Pc$ denote its image under the embedding $\PrG\to\ParG$. Then, $\pi_1(B\Pc)\cong U(P)$.
\end{cor}

\begin{dem}
Through the embedding of Proposition \ref{prop pregroup is partial group}, the universal property satisfied by $U(P)$ (Proposition \ref{prop:universalgroup}) becomes exactly the one satisfied by $\pi_1(B\Pc)$ (Proposition \ref{prop:fundamental group}).  
\end{dem}

\begin{rem}
At first glance, one may wonder whether this construction of a partial group from a pregroup could be reciprocal: if we restrict the domain of a given partial group to words of length two, do we actually get a pregroup? A second look at the axioms for pregroups might convince you that the answer is ``no'' in general: the Axioms \ref{cond:Pre3} and \ref{cond:Pre4} allow to deduce that some words belong to the domain without already knowing that a bigger word containing them belongs to it too, and this kind of property does not seem to appear in the definition of a partial group. 
\end{rem}

Now, we give two examples of localities (cf. Definition \ref{def locality}). The first example is actually a pregroup, but the second one is not, so that it provides an example of a partial group which is not a pregroup. 

\begin{ex}
Let $G=GL_3(\Fb_2)$ and let $S \leq G$ be the subgroup of upper-triangular matrices with diagonal coefficients equal to 1. Then $S$ is a 2-subgroup of $G$ isomorphic to the dihedral group of order 8. As such, it contains three subgroups of order $4$: two of them, that we will denote by $V$ and $V'$, are isomorphic to $C_2\times C_2$, and the last one is cyclic and denoted by $C$. Set $\Delta=\{C,V,V',S\}$ and $\Lc=\{g\in G\mid \exists Q\in \Delta, \ Q^g\in \Delta\}$. Then one can check that, taking $\Db_\Delta$ to be the domain (cf. Definition \ref{def objective}), $(\Lc,\Delta,S)$ is a locality. Also, the subgroups $V$ and $V'$ are not conjugate in $G$ and thus not in $\Lc$, and both are clearly not conjugate to $C$. Moreover, $N_G(V)\cap N_G(V')=S=N_G(C)$. Hence one gets that $\Db_\Delta=\Wb(N_G(V))\cup\Wb(N_G(V'))$, so, by Example \ref{ex:A*_CB revisted}, $\Lc$ is obtained as in Example \ref{ex:A*_CB}, i.e. $\Lc = N_G(V) \cup_S N_G(V')$.
\end{ex} 

\begin{ex}
Let $\Tc$ be one of the linking systems constructed in \cite{COS} such that the fundamental group of the geometric realisation of $\Tc$ is trivial. The associated locality $\Lc$ (through the correspondence highlighted by Chermak in \cite[Appendix A]{Ch0}) also satisfies that $\pi_1(B\Lc)$ is trivial thanks to \cite[Theorem A.5]{Go}. Hence if $\Lc$ were a pregroup, then by Corollary \ref{cor:universal=pi1}, $U(\Lc)\cong \pi_1(B\Lc)$ would be trivial, and thus $\Lc$ would be too. However, $\Lc$ is clearly non trivial, so $\Lc$ is not a pregroup.
\end{ex}

\section{Graphs of groups}

\label{Section3}

\subsection{Graphs of groups and their fundamental group}

\begin{defi}
A \emph{graph} $Y=(V,E)$ is the data of:
\begin{enumerate}[$(i)$]
\item a set $V$ of \emph{vertices},
\item a set $E$ of \emph{edges},
\item two maps $\iota : E \to V$ and $\tau : E \to V$ mapping each edge $e$ to its \emph{initial vertex} and \emph{terminal vertex} respectively,
\item a fixed point-free involution of the edges, denoted $e \mapsto \overline{e}$, such that for all $e \in E$, $\iota(e) = \tau(\overline{e})$.
\end{enumerate}
A graph is said to be \emph{connected} if for every $x,y \in V$ there exist edges $e_0,e_1,\cdots ,e_n$ such that $\iota(e_0) = x$, $\tau(e_n) = y$ and for all $i \in \{ 0,1, \cdots , n-1\}$, $\tau(e_i) = \iota(e_{i+1})$.
\end{defi}

\begin{defi}
A \emph{graph of groups} $(G,Y)$ consists of a connected graph $Y=(V,E)$, a group $G_v$ for every vertex $v \in V$, and a group $G_e$ for every edge $e \in E$, together with a monomorphism $G_e \to G_{\tau(e)}$ denoted by $g \mapsto g^e$, such that $G_{\overline{e}}=G_e$. The $G_v$ and $G_e$ are called \emph{vertex groups} and \emph{edge groups} respectively, and the maps $g \mapsto g^e$ and $g \mapsto g^{\overline{e}}$ are called the \emph{edge maps}. Thus there is one edge group and two edge maps by \emph{geometric edge}, i.e. by pair $(e,\overline{e})$.
\end{defi}

To every graph of group, we can associate a particular group called its fundamental group.

\begin{defi}\label{def fundamental group}
Let $(G,Y)$ be a graph of groups, with $Y=(V,E)$, and let $T=(V,E')$ be a maximal tree in $Y$. The \emph{fundamental group} of $(G,Y)$ is the group generated by the vertex groups $G_v$ for $v \in V$ and by the edges $e \in E$, subject to the following relations: \begin{enumerate}[$\bullet$]
\item $e^{-1} = \overline{e}$ for each $e \in E$,
\item $e a^e e^{-1} = a^{\overline{e}}$ for each $e \in E$ and each $a \in G_e$,
\item $e = 1$ for each $e \in E'$.
\end{enumerate}
\end{defi}

We speak about ``the'' fundamental group of $(G,Y)$ and not about the fundamental group of $(G,Y)$ with respect to the tree $T$ because it is in fact independent of the choice of $T$. This result is a direct consequence of \cite[I, Proposition 20]{Se}.
\begin{prop}
The fundamental group of a graph of groups $(G,Y)$, as defined above, does not depend on the choice of the maximal tree in $Y$.
\end{prop}

Several classical constructions in combinatorial group theory arise as fundamental groups of particular graphs of groups. Let us mention two of them.

\begin{exs} \label{examples graphs of groups}
\textbf{(1)} \quad Let $A$ and $B$ be two groups, and consider a third group $C$ with two monomorphisms $\varphi \colon C \to A$ and $\psi \colon C \to B$. This forms a graph of groups whose underlying graph has two vertices with one geometric edge between them, the vertex groups being $A$ and $B$, and the edge group being $C$ with edge maps $\varphi$ and $\psi$ (see Figure \ref{fig:A*B}). Then the fundamental group of this graph of groups is the \emph{amalgamated product of $A$ and $B$ over $C$} (or the \emph{free product of $A$ and $B$ amalgamating $C$}), denoted $A*_C B$. It equals $(A*B)/N$, where $N$ is the normal subgroup of $A*B$ generated by all relations of the form $\varphi(c)\psi(c)^{-1}$ for $c \in C$.

\begin{figure}[!h]
\centering
\begin{tikzpicture}[node distance=2.5cm]

\tikzstyle{group}=[draw,circle,minimum size=8pt,inner sep=0pt]
 \node[group] (A)  [label=left:$A$] {};
 \node[group] (B)  [right of=A,label=right:$B$]   {};
 
 \draw[-,color=black] (A) edge [out=0,in=180] node[above] {$C$} (B);
\end{tikzpicture}
\caption{A graph of groups whose fundamental group is $A*_C B$.\label{fig:A*B}}
\end{figure}
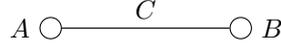

\medskip

\textbf{(2)} \quad Let $G$ be a group, $H$ a subgroup of $G$, and let $\alpha \colon H \to G$ be a group monomorphism. We construct a graph of groups by taking only one vertex, with vertex group $G$, and one geometric edge, with edge group $H$ (see Figure \ref{fig:HNN}). We take $\alpha$ and the inclusion $H \hookrightarrow G$ as edge maps. The fundamental group of this graph of groups is the \emph{HNN extension of $G$ by $\alpha$}, denoted $G*_\alpha$. It is the group $(G*\Zb)/N$ where $N$ is the normal subgroup generated by all relations of the form $tht^{-1}\alpha(h)^{-1}$ for $h \in H$, $t$ being a fixed generator of $\Zb$ called the \emph{stable letter} of the HNN extension. 

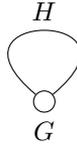
\begin{figure}[!h]
\centering
\begin{tikzpicture}[node distance=2.5cm]
\tikzstyle{group}=[draw,circle,minimum size=8pt,inner sep=0pt]
 \node[group] (G) [label=below:$G$] {};
 
 \draw[-,color=black] (G) edge [out=40,in=140,looseness=20] node[above] {$H$} (G);
\end{tikzpicture}
\caption{A graph of groups whose fundamental group is $G*_\alpha$. \label{fig:HNN}}
\end{figure}
\end{exs}

It is well-known (cf. \cite[I, Theorem 8]{Se}) that any finite subgroup of the amalgamated product $A*_C B$ is conjugate to a subgroup of $A$ or $B$. It is in fact a general property of fundamental groups of graphs of groups.

\begin{prop}[{\cite[Corollary 28]{LS}}] \label{prop finite subgroups graph of groups}
Every finite subgroup of the fundamental group of a graph of groups is conjugate to a subgroup of a vertex group.
\end{prop}

\begin{rem}
The article \cite{LS}, which is our reference for the above proposition, works with a topological definition for the fundamental group of a graph of groups, which is apparently different from ours. However both definitions lead to isomorphic groups, as can be seen in \cite{LiSe}. In a nutshell: the topological approach defines the fundamental group of a given graph of groups $(G,Y)$ as the $\pi_1$ of the homotopy colimit of the graph of spaces obtained by application of the classifying-space functor $B(-)$ on $(G,Y)$. Proposition 3.2 in \cite{LiSe} then tells us that this homotopy colimit is homeomorphic to the classifying space of the fundamental group as defined in Definition \ref{def fundamental group}. Taking the $\pi_1$ on both sides gives us the desired isomorphism.
\end{rem}

\subsection{Graphs of groups and universal groups of pregroups}

In \cite{Ri}, Rimlinger makes an extensive study of the relations between pregroups and graphs of groups. In fact, he proved that universal groups of pregroups and fundamental groups of graphs of groups are closely related constructions: the universal group of a pregroup $P$ is, under one condition on $P$ (that of being of finite height - which includes finite pregroups and a lot more), the fundamental group of a particular graph of groups constructed from $P$. Reciprocally, the fundamental group of a graph of groups $\Gc$ is, under some conditions on $\Gc$, the universal group of a particular pregroup constructed from $\Gc$. For our purpose, we only need these results for finite pregroups and finite graphs of finite groups, so we state them in this particular case.

\begin{thm}[{cf. \cite[Theorem A]{Ri}}] \label{thm rimlinger 1}
Let $P$ be a finite pregroup. We denote $U(P)$ its universal group. Then there exists a finite graph of finite groups $\Gc$ whose edge groups are subgroups of $P$ and whose fundamental group is isomorphic to $U(P)$.
\end{thm}

From Theorem \ref{thm rimlinger 1} and Proposition \ref{prop finite subgroups graph of groups}, we deduce the following result.
\begin{cor} \label{cor finite subgroups pregroups}
Let $P$ be a finite pregroup and $U(P)$ its universal group. Then every finite subgroup of $U(P)$ is conjugate to a subgroup of $P$.
\end{cor}

For the reciprocal of Theorem \ref{thm rimlinger 1}, apart from being finite, we need another condition on the graph of groups.

\begin{defi}
Let $(G,Y)$ be a graph of groups. We say that $(G,Y)$ is \emph{proper} if none of its edge maps are surjective.
\end{defi}

\begin{thm}[{cf. \cite[Theorem B]{Ri}}] \label{thm rimlinger 2}
Let $(G,Y)$ be a finite graph of finite groups which is proper. Then there exists a finite pregroup $Q$ whose universal group is isomorphic to the fundamental group of $(G,Y)$.
\end{thm}

\begin{rem}
As the statement of Theorem B in \cite{Ri} is fairly indigestible, the careful reader is entitled to ask how we can see that the pregroup is finite whenever $(G,Y)$ is a finite graph of finite pregroups. Let us give some elements of exegesis of Rimlinger's proof to reassure our reader. Here we refer only to \cite{Ri} and stick to its notations (exept that we denoted our graph of groups $(G,Y)$ instead of $(\bm{H},Y)$). The pregroup $Q$ constructed for Theorem B is defined in Definition 7.13 as the preimage in $F(G,Y)$ of a finite set of paths in $Y$ (in bijection with the disjoint union of the set of vertices and the set of edges outside a maximal subtree) under a certain map also denoted $Y$. For our explanation, one only needs to know that $F(G,Y)$ is the universal group of a pregroup $P$ (Theorem 7.7) which is the quotient under a certain equivalence relation of the pregroup $P'$ defined at the beginning of Chapter 7. In the case where $(G,Y)$ is a finite graph of finite groups, the definition of the pregroup $P'$ clearly implies that it is finite (hence $P$ is finite too). Moreover $Y$ is a map defined on $P'$ which is compatible with the quotient $P' \to P$, and extends to words in $P'$ via $(x_1,\cdots,x_n) \mapsto \big(Y(x_1),\cdots,Y(x_n) \big)$. By Lemma 7.9, the map $Y$ is well-defined on $F(G,Y)$ considering its value on $P$-reduced representative $P$-words. Now, given a fixed path $(y_1,\cdots,y_n)$, for each $i \in \{1,\cdots,n\}$ there exists only a finite number of $x_i \in P$ such that $Y(x_i)=y_i$ because $P$ is finite. Thus there exists only a finite number of ($P$-reduced) $P$-words $(x_1, \cdots , x_n)$ such that $Y(x_1,\cdots,x_n) = (y_1 , \cdots , y_n)$. As $Q$ is defined to be the preimage of a finite set of paths under $Y$, we can conclude that $Q$ is finite.
\end{rem}

\section{Fusion systems and realisability}

\label{Section4}

\subsection{Definitions and examples}

Let $p$ be a fixed prime number. A fusion system over a finite $p$-group $S$ is a way to abstract the action of a group $G$ containing $S$ on the subgroups of $S$ by conjugation. Given a group $G$ and an element $g\in G$, we will denote by $c_g$ the homomorphism $x\in G\mapsto g^{-1}xg\in G$ (this is consistent with Notation \ref{notation conjugaison}). Our convention for the composition of two maps $f \colon X \to Y$ and $f' \colon Y \to Z$ is to denote it by $f' \circ f$, so that $c_{g_1} \circ c_{g_2}$ equals $c_{g_2g_1}$. For $H,K$ two subgroups of $G$, $\Hom_G(H,K)$ will denote the set of all group homomorphisms $c_g$, for $g\in G$ such that $c_g(H) \leq K$. Finally, $\Syl_p(G)$ will denote the collection of all Sylow $p$-subgroups of $G$.

\begin{defi}\label{defF}
Let $S$ be a finite $p$-group. A \emph{fusion system} over $S$ is a small category $\Fc$ whose object set $\Ob(\Fc)$ is the set of all subgroups of $S$ and whose morphism sets $\Mor_\Fc(P,Q)$, for $P,Q \leq S$, satisfy the following two properties:
 \begin{enumerate}[label = (F\arabic*)]
  \item\label{F1} $\Hom_S(P,Q)\subseteq \Mor_\Fc(P,Q)\subseteq \Inj(P,Q)$;
  \item\label{F2} each $\varphi\in\Mor_\Fc(P,Q)$ is the composite of an $\Fc$-isomorphism followed by an inclusion.
 \end{enumerate}
\end{defi}

The composition law in a fusion system is given by composition of homomorphisms. We usually write $\Hom_\Fc(P,Q):=\Mor_\Fc(P,Q)$ to emphasise the fact that the morphisms in $\Fc$ are group homomorphisms.

\begin{rems}
\begin{enumerate}[$(1)$]
\item Over a fixed $p$-group $S$, there is a minimal fusion system. Its morphism sets are the $\Hom_S(P,Q)$, for $P,Q \leq S$. It is called the \emph{inner fusion system} of $S$ and denoted $\F_S(S)$. There is also a maximal fusion system over $S$, with morphism sets equal to $\Inj(P,Q)$ for $P,Q \leq S$.
\item The \emph{intersection} of two (or more) fusion systems over the same $p$-group $S$ is obtained by taking the intersection of the morphism sets for each fusion system. This forms again a fusion system. Thus it makes sense to talk about the fusion system over $S$ \emph{generated} by a certain family of injective group homomorphisms between subgroups of $S$.
\end{enumerate}
\end{rems}

The typical example of a fusion system is the fusion system of a finite group $G$ over one of its Sylow $p$-subgroups, but we can define more generally the fusion system of any group over one of its $p$-subgroups.

\begin{ex}\label{ex:FusGroup}
Let $S$ be a finite $p$-group, and let $G$ be a group containing $S$. The \emph{fusion system of $G$ over $S$} is the category $\Fc_S(G)$ where $\Ob(\Fc_S(G))$ is the set of all subgroups of $S$ and $\Mor_{\Fc_S(G)}(P,Q)=\Hom_G(P,Q)$ for all $P,Q\leq S$. One can easily check that the category $\Fc_S(G)$ defines a fusion system over $S$.
\end{ex}

In fact we can generalise even further, considering the fusion system of a partial group over one of its $p$-subgroups. We use Notation \ref{notation conjugaison} in the definition.

\begin{defi} \label{def fusion partial group}
Let $\Mc$ be a partial group and let $S$ be a $p$-subgroup of $\Mc$. We define the \emph{fusion system of $\Mc$ over $S$}, denoted $\F_S(\Mc)$, to be the fusion system over $S$ generated by conjugation maps $c_g \colon Q \to Q^g$, whenever $g \in \Mc$ is such that $Q \subseteq \Db(g)$, $Q^g$ is a subgroup of $S$ and $c_g \colon Q \to Q^g$ is a group homomorphism. For $\F_S(\Mc)$ to be well-defined, one only needs to check that $c_g$ is injective, which is already known from Lemma \ref{lem partial conjugaison}.
\end{defi}

As a particular case, we recover the fusion system of a locality, which was introduced by Chermak in \cite{Ch0}. Notice that if $\Mc = \Lc$ is a locality associated to $S$, whenever $c_g \colon Q \to S$ is defined, then $Q^g$ is a subgroup of $S$ and $c_g$ is a group homomorphism (cf. Proposition 2.6 in \cite{Ch1}). Hence, morphisms in $\F_S(\Lc)$ are then just compositions of restrictions of conjugation maps $c_g$ (for $g \in \Lc$) between subgroups of $S$. This does not seem to be the case in general partial groups.

The fusion system of a locality, as well as the fusion system of a finite group over one of its Sylow $p$-subgroups, belong to the important family of \emph{saturated fusion systems}. This notion will not be discussed here (see for example \cite{AKO}), but the idea is that $\F$ ``behaves'' like $\F_S(G)$ when $S$ is a Sylow $p$-subgroup of a finite group $G$. In the literature, when a saturated fusion system is isomorphic to a fusion system of the form $\F_S(G)$ for $G$ a finite group and $S \in \Syl_p(G)$, the fusion system is said to be \emph{realisable} (it is called \emph{exotic} otherwise). In the following we enlarge this notion, discussing about the ``realisability'' of general fusion systems in some subclass of the class of partial groups. First, we need to define a notion of Sylow $p$-subgroups for partial groups.

\begin{defi}
Let $\Mc$ be a partial group and let $S$ be a $p$-subgroup of $\Mc$. We say that $S$ is a \emph{Sylow $p$-subgroup} of $\Mc$ if for every $p$-subgroup $P$ of $\Mc$ there exists a sequence $(g_1, \cdots , g_r)$ of elements of $\Mc$ and a sequence $(P_0,\cdots , P_r)$ of $p$-subgroups of $\Mc$ such that: \begin{enumerate}[$\bullet$]
\item $P_0 = P$;
\item for each $i \in \{1 , \cdots , r\}$, $P_i \subseteq \Db(g_i)$, $c_{g_i} \colon P_i \to {P_i}^{g_i}$ is a group homomorphism, and ${P_i}^{g_i} = P_{i+1}$;
\item $P_r \leq S$.
\end{enumerate}
\end{defi}

\begin{rems} \label{rem p-sylow def}
\begin{enumerate}[$(1)$]
\item A Sylow $p$-subgroup of $\Mc$, if it exists, is a $p$-subgroup of maximal order. However it could happen that $S$ is a $p$-subgroup of maximal order in $\Mc$ but not a Sylow $p$-subgroup. For example, if we consider the pregroup $P=A \cup_C B$ as in Example \ref{ex:A*_CB}, with $A = C_2$, $B = C_4$ and $C=\triv$, then $B$ is a $2$-subgroup of $P$ of maximal order but the $2$-subgroup $A$ is not conjugate to any subgroup of $B$ because the only conjugation maps in $P$ defined on $A$ are conjugation by elements in $A$.
\item If $\Mc$ is a group, this definition is equivalent to asking that any $p$-subgroup of $\Mc$ is conjugate in $\Mc$ to some subgroup of $S$. In particular when $\Mc$ is a finite group, we recover the classical definition of Sylow $p$-subgroup, by the Sylow theorems.
\item If $\Mc=\Lc$ is a locality, Andrew Chermak (see \cite[Definition 2.16]{Ch1}) also gave a definition of a Sylow $p$-subgroup $S$, which asks for the existence of a set $\Delta$ of subgroups of $S$ such that $(\Lc,\Delta,S)$ is again a locality. Our definition is a priori broader, as when $(\Lc,\Delta,S)$ is a locality, every $p$-subgroup of $\Lc$ is conjugate to a subgroup of $S$ by \cite[Proposition 2.11.(c)]{Ch1}.
\end{enumerate}
\end{rems}

\begin{defi}
Let $\Cc$ be a subclass of the class of partial groups. Given a fusion system $\F$ over a finite $p$-group $S$, we say that $\F$ is \emph{realisable in $\Cc$} if there exists an objet $X$ in $\Cc$, containing $S$ as a Sylow $p$-subgroup, such that $\F$ is isomorphic to $\F_S(X)$. We say that $\F$ is \emph{weakly realisable in $\Cc$} if there exists an objet $X$ in $\Cc$ containing $S$ such that $\F$ is isomorphic to $\F_S(X)$.
\end{defi}

Of course, this definition will only be of interest for some particular classes $\Cc$. As examples, we restate several results from the literature in these terms.

\begin{exs}
\begin{enumerate}[$(a)$]
\item Let $\Cc$ be the class of finite groups. If we consider the property of being realisable in $\Cc$ for saturated fusion systems, we recover the classical use of the terminology. Saturated fusion systems which are not realisable in $\Cc$ are the so-called \emph{exotic} fusion systems. Moreover, in \cite{Pa}, Sejong Park proved that any fusion system is weakly realisable in $\Cc$. 
\item Let $\Cc$ be the class of localities. In \cite[Main theorem]{Ch0}, Andrew Chermak proved that any saturated fusion system is realisable in $\Cc$. In this case, the Sylow $p$-subgroup of the locality is even a Sylow $p$-subgroup in the sense of Chermak (see $(3)$ in Remarks \ref{rem p-sylow def} above).
\item Let $\Cc$ be the class of (not necessarily finite) groups. Leary \& Stancu, in \cite{LS}, proved that any fusion system is realisable in $\Cc$. Independently at the same time, Robinson (in \cite{Ro}), using a different construction, also proved that a large class of fusion systems (including saturated fusion systems) is realisable in $\Cc$.
\item Let $\Cc = \PrG$ be the class of finite pregroups. In the following section we prove that any fusion system is realisable in $\Cc$.
\end{enumerate}
\end{exs}

We detail here a slightly adapted version of the result of Leary and Stancu mentioned in Example $(c)$ above, in view of a further use.

\begin{thm}[cf. {\cite[Theorem 2]{LS}}] \label{thm Leary-Stancu}
Let $\F$ be a fusion system over a finite $p$-group $S$. Assume that $\F$ is generated by $\Phi = \{ \Phi_1, \cdots , \Phi_r \}$, where each $\Phi_i$ is an injective group homomorphism $P_i \to Q_i$ between subgroups of $S$. Let $U$ be any finite group whose order is prime to $p$. Let $G$ be the iterated HNN-extension $\big(\cdots \big((S\times U)*_{\Phi_1} \big)*_{\Phi_2} \cdots \big)*_{\Phi_r}$. Then $S$ embeds as a Sylow $p$-subgroup of $G$, and $\F_S(G) = \F$.
\end{thm}

\begin{figure}[!h]
\centering
\begin{tikzpicture}[node distance=1cm]

\tikzstyle{group}=[draw,circle,minimum size=8pt,inner sep=0pt]
 \node[group] (S)  [label=below:$S \times U$] {};
 \node             [above right=of S] {$\ddots$};
 
 \draw[-,color=black] (S) edge [out=195,in=165,looseness=100] node[label=180:$P_1$] {} (S);
 \draw[-,color=black] (S) edge [out=165,in=135,looseness=100] node[label=150:$P_2$] {} (S);
 \draw[-,color=black] (S) edge [out=135,in=105,looseness=100] node[label=120:$P_3$] {} (S);
 \draw[-,color=black] (S) edge [out=105,in=75,looseness=100] node[label=90:$P_4$] {} (S);
 \draw[-,color=black] (S) edge [out=15,in=345,looseness=100] node[label=0:$P_r$] {} (S);
\end{tikzpicture}
\caption{A graph of groups whose fundamental group is the group $G$ in Theorem \ref{thm Leary-Stancu}.}
\label{figure graphe Leary-Stancu}
\end{figure}
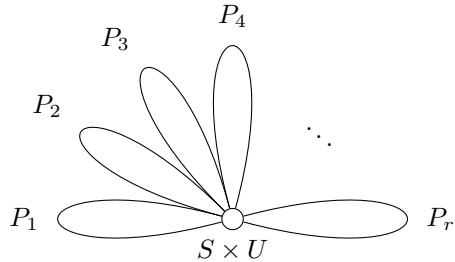

The group $G$ in Theorem \ref{thm Leary-Stancu} above is the fundamental group of the graph of groups in Figure \ref{figure graphe Leary-Stancu}, where for each $i \in \{1, \cdots ,r\}$, the two edge maps $P_i \to S \times U$ are $\Phi_i$ and the inclusion $P_i \hookrightarrow S$, both post-composed with the inclusion $S \hookrightarrow S \times U$.

\begin{dem}[sketch]
Most of the proof of \cite[Theorem 2]{LS} remains valid \textit{mutatatis mutandi} if we replace $S$ with $S \times U$. Let us just mention the two small arguments that need to be added. First, any finite $p$-subgroup of $G$ being conjugate to a subgroup of $S \times U$, it is in fact conjugate to a subgroup of $S$ because $|U|$ is prime to $p$, so $\F_S(G)$ remains well-defined. Secondly, when one takes an element $g \in G$ which conjugates two subgroups of $S$, the proof of \cite[Theorem 2]{LS} tells us that $c_g$ is a composition of morphisms of the form $\Phi_i$ or $\Phi_i^{-1}$ for some $i \in \{1, \cdots , r\}$ and of conjugation morphisms $c_t$ for some $t \in S \times U$, defined between subgroups of $S \times U$. But in fact, since we start with a subgroup of $S$, we must arrive on a $p$-subgroup of $S \times U$ (i.e. a subgroup of $S$) at each composition step. Moreover, conjugation by an element of $S \times U$ between subgroups of $S$ is nothing but conjugation by an element of $S$, because $U$ commutes with $S$. Thus, $c_g$ can be expressed as a composition of morphisms of the form $\Phi_i$ or $\Phi_i^{-1}$ for some $i \in \{1, \cdots , r\}$ and of conjugation morphisms $c_s$ for some $s \in S$, all defined between subgroups of $S$. In other words, $c_g$ belongs to $\F$.
\end{dem}

\medskip

Finally, we will use the notion of morphism of fusion systems.

\begin{defi}
Let $\F$ and $\F'$ be two fusion systems over finite $p$-groups $S$ and $S'$ respectively. A \emph{morphism of fusion systems} from $\F$ to $\F'$ is a group homomorphism $\alpha_0 \colon S \to S'$ together with a covariant functor $\alpha \colon \F \to \F'$, so that $\alpha(P) = \alpha_0(P)$ and $\big(\alpha (\varphi) \circ \alpha_0 \big) (u) = (\alpha_0 \circ \varphi) (u)$ for any $P \leq S$, any $\varphi \colon P \to S$ in $\F$ and any $u \in P$.
\end{defi}

Fusion systems, together with morphisms of fusion systems and the usual composition of both group homomorphism and functors, form a category. In particular, we are usually interested in fusion systems up to isomorphism in this category.

\begin{lem} \label{lem isomorphisme systèmes de fusion}
Let $G$ be a group, let $S$ be a $p$-subgroup of $G$. For any $g \in G$, the fusion systems $\F_S(G)$ and $\F_{S^g}(G)$ are isomorphic.
\end{lem}

\begin{dem}
Consider the group homomorphism $c_g \colon S \to S^g$ and the functor ${c_g}^* \colon \F_S(G) \to \F_{S^g}(G)$ defined on subgroups of $S$ by $P \mapsto P^g$. The action of ${c_g}^*$ on a morphism $\varphi \colon P \to Q$ in $\F_S(G)$ is given by $$ {c_g}^*(\varphi) \ = \ c_g \circ \varphi \circ c_g^{-1} \colon P^g \longrightarrow Q^g \ .$$ This clearly defines a morphism of fusion systems between $\F_S(G)$ and $\F_{S^g}(G)$. Moreover, $c_{g^{-1}} \colon S^g \to S$ is the inverse of $c_g$ and one can check that ${c_{g^{-1}}}^*$ is an inverse for ${c_g}^*$, so this morphism of fusion systems is an isomorphism.
\end{dem}

\subsection{Realisability of fusion systems in finite pregroups}

As we proved in Proposition \ref{prop pregroup is partial group}, pregroups are particular instances of partial groups. Thus, we can construct the fusion system of a pregroup over one of its (Sylow) $p$-subgroups and ask if any fusion system can be obtained in this way, i.e. if any fusion system is realisable in the class of finite pregroups.

First, considering a pregroup $(P,D)$ as a partial group through the embedding of Proposition \ref{prop pregroup is partial group} and using Notation \ref{notation conjugaison}, it is straightforward to check that if $g \in P$ we have \begin{align*}
\Db(g) \ &= \ \{ x \in P \ \vert \ (g^{-1},x),(x,g) \in D \textup{ and } (g^{-1}x,g) \in D \} \\
&= \ \{ x \in P \ \vert \ (g^{-1},x),(x,g) \in D \textup{ and } (g^{-1},xg) \in D \} \quad \textup{thanks to Axiom \ref{cond:Pre3}} \\
&= \ \{ x \in P \ \vert \ (g^{-1},x,g) \textup{ associates}\} \ . 
\end{align*}

The following lemma shows that conjugation maps are always group homomorphisms in the context of pregroups.

\begin{lem} \label{lem conjugation morphisms}
Let $(P,D)$ be a pregroup. Let $g \in P$ and let $Q$ be a subgroup of $P$ such that $Q \subseteq \Db(g)$. We denote $Q^g$ the image of $Q$ by $c_g$ as in Notation \ref{notation conjugaison}. Then $Q^g$ is a subgroup of $P$, contained in $\Db(g^{-1})$, and $c_g \colon Q \to Q^g$ is a group isomorphism whose inverse is $c_{g^{-1}} \colon Q^g \to Q$.
\end{lem}

\begin{dem}
It is enough to prove $1_P \in Q^g$, $Q^g \times Q^g \subseteq D$ and that $Q^g$ is a group for the operations iduced by those of $P$ to conclude that it is a subgroup of $P$. As $1_P \in Q$ and $g^{-1}1_Pg = 1_P$, it is clear that $1_P \in Q^g$. Moreover, the fact that $Q^g$ is stable under $(-)^{-1}$ comes from items \ref{lem item involution} and \ref{lem item inverse} in Lemma \ref{lem propriétés prégroupe}. 

Now let $x_1,x_2 \in Q$, hence $(x_1,x_2) \in D$ and $x_1x_2 \in Q$. We want to prove that $(x_1^g,x_2^g) \in D$ and $x_1^g x_2^g \in Q^g$. On the one hand we have $(g^{-1}x_1,g),(g^{-1},x_2g) \in D$, so by item \ref{lem item intercalation} in Lemma \ref{lem propriétés prégroupe} we have the following equivalence: $$ (x_1^g,x_2^g) \in D \iff (g^{-1}x_1,x_2g) \in D $$
On the other hand, using $x_1x_2 \in Q$ and two applications of Axiom \ref{cond:Pre3}, we know that $$ (g^{-1},(x_1x_2)g) \in D \iff (g^{-1},x_1(x_2g)) \in D \iff (g^{-1}x_1,x_2g) \in D \ .$$ Thus $(x_1^g,x_2^g) \in D$, and repeated use of Axiom \ref{cond:Pre3} gives $$ (x_1^g)(x_2^g) \ = \ g^{-1} (x_1x_2) g \ \in Q^g.$$
The above equality also tells us that $c_g \colon Q \to Q^g$ is a group homomorphism. The rest of the proof is just Lemma \ref{lem partial conjugaison}.
\end{dem}

\begin{rem}
If $P$ be a pregroup and $S$ is a $p$-subgroup of $P$, then the morphisms of the fusion system $\F_S(P)$, as defined in Definition \ref{def fusion partial group}, are precisely the maps $Q \to R$ between subgroups of $S$ which are composition of restrictions of conjugation maps $c_g$ for $g \in P$.
\end{rem}

Now we can prove that the fusion system of a finite pregroup $P$ and that of its universal group $U(P)$, over the same Sylow $p$-subgroup $S$, coincide.

\begin{thm} \label{thm système prégroupe}
Let $P$ be a pregroup and $S$ be a finite $p$-group. Then $S$ embeds as a Sylow $p$-subgroup of $P$ if and only if it embeds as a Sylow $p$-subgroup of $U(P)$, and in this case we have $\F_S \big(U(P) \big)=\F_S(P)$. 
\end{thm}

\begin{dem}
We identify $P$ with its image through the canonical embedding $P \hookrightarrow U(P)$. It is clear that any Sylow $p$-subgroup $S$ of $P$ is a $p$-subgroup of $U(P)$. According to Corollary \ref{cor finite subgroups pregroups}, any $p$-subgroup $Q$ of $U(P)$ is conjugate to a $p$-subgroup $R$ of $P$ and hence there exists a sequence $Q=Q_0,R=Q_1, Q_2, \cdots , Q_r$ of subgroups of $U(P)$ such that each $Q_{i+1}$ is conjugate (in $U(P)$) to $Q_i$ and $Q_r \leq S$. So $S$ is in fact a Sylow $p$-subgroup of $U(P)$. 

Reciprocally, if $S$ is a Sylow $p$-subgroup of $U(P)$, then by Corollary \ref{cor finite subgroups pregroups} it embeds in $P$ via a conjugation morphism $c_g$, with $g \in U(P)$. Then any $p$-subgroup of $P$ being a $p$-subgroup of $U(P)$, it is conjugate to a subgroup of $S$, and hence conjugate in $U(P)$ to a subgroup of $S^g$. We then deduce what we need from Proposition \ref{prop conjugate sgs in P}, so that $S^g$ is a Sylow $p$-subgroup of $P$.

By Lemma \ref{lem isomorphisme systèmes de fusion}, the fusion systems $\F_S\big(U(P)\big)$ and $\F_{S^g}\big(U(P) \big)$ are isomorphic for all $g \in G$. We can thus assume that $S$ is a Sylow $p$-subgroup of $P$. The fusion system $\F_S \big(U(P) \big)$ clearly contains all the morphisms in $\F_S(P)$. Since $\F_S \big(U(P) \big)$ is generated by conjugation maps $c_g \colon Q \to R$ with $g \in U(P)$ and $Q,R \leq S$, it is enough to prove that any such morphism belongs to $\F_S(P)$. We can assume that $R = c_g(Q)$. Now Proposition \ref{prop conjugate sgs in P} precisely tells us that $c_g$ is equal to some composition of conjugation maps $c_{g_i}$ defined between subgroups of $S$, with $g_i \in P$ for each $i$. In other words, $c_g$ belongs to $\F_S(P)$, which concludes the proof.
\end{dem}

Combining this theorem with the result of Leary \& Stancu and the second theorem of Rimlinger, we can prove that every fusion system is realisable in the class of finite pregroups.

\begin{cor}\label{cor:realisability pregroups}
Every fusion system $\F$ over a finite $p$-group $S$ is the fusion system of a finite pregroup containing $S$ as a Sylow $p$-subgroup. 
\end{cor}

\begin{dem}
Suppose that $\F$ is generated by $\{ \Phi_1, \cdots , \Phi_r \}$, where $\Phi_i$ is a morphism defined on a subgroup $P_i$ of $S$ for each $i \in \{ 1 , \cdots , r \}$. Let $q$ be a prime number different from $p$ and let $C_q$ denote the cyclic subgroup of order $q$ (but we could take any non-trivial finite group whose order is prime to $p$ instead). Set $G := \big(\cdots ((S\times C_q)*_{\Phi_1})*_{\Phi_2} \cdots \big)*_{\Phi_r}$. Then $S$ embeds in $G$, and $G$ is the fundamental group of the graph of groups in Figure \ref{figure graphe Leary-Stancu}.

This finite graph of finite groups is proper (that was the whole purpose of adding the $C_q$ part to $S$), so by Theorem \ref{thm rimlinger 2} there exists a finite pregroup $P$ whose universal group is isomorphic to $G$. Moreover, according to Theorem \ref{thm Leary-Stancu}, $G$ contains $S$ as a Sylow $p$-subgroup and the fusion system of $G$ over $S$ is nothing but $\F$. 

Finally, by Theorem \ref{thm système prégroupe}, we get that $S$ embeds as a Sylow $p$-subgroup of $P$, and $$ \F_S(G) \ = \ \F_S \big(U(P) \big) \ = \ \F_S(P) \ . $$
\end{dem}

\section{Examples of pregroups realising fusion systems}\label{sec:Examples}

In this section, we detail two constructions of pregroups realising fusion systems, according to Theorem \ref{thm système prégroupe}. The first construction leads to a pregroup whose universal group is the ``Leary-Stancu group'' given in Theorem \ref{thm Leary-Stancu}, thus providing a more direct way to prove Corollary \ref{cor:realisability pregroups} (without refering to Theorem \ref{thm rimlinger 2}). The second construction similarly leads to a pregroup whose universal group is the ``Robinson group'' used in \cite[Theorem 2]{Ro}.

\subsection{A pregroup for the Leary-Stancu group}

Let $\F$ be any fusion system over a finite pregroup $S$. Assume that $\F$ is generated by a certain family of morphisms $\{ \phi_1, \cdots , \phi_r \}$, where each $\phi_i$ is a group isomorphism $P_i \to Q_i$, with $P_i,Q_i \leq S$. We associate a symbol $t_i$ to each $\phi_i$, $i \in \{1, \cdots ,r \}$, and we take the free product $F$ of $S$ with the free group generated by $\{ t_1, \cdots , t_r\}$. The Leary-Stancu group $G$ is the quotient of $F$ by the normal closure of the elements ${t_i}^{-1}ut_i\phi_i(u)^{-1}$, where $i \in \{1,\cdots ,r\}$ and $u \in P_i$. The elements $t_i \in G$ will sometimes be referred to as the \emph{stable letters} of $G$.

\subsubsection{Constructing the pregroup}

Informally, we consider the subset $P$ of $G$ formed by the elements $s \in S \leq G$ and all the elements of the form $a t_i a'$ or $b t_i^{-1} b'$, with $a,a',b,b' \in S$ and $i \in \{1, \cdots ,r\}$. A pair $(x,y) \in P\times P$ belongs to the domain $D$ if and only if the product $xy$ in $G$ belongs to $P$, and we define the inverse and multiplication in $P$ as it is in $G$.

More formally, for each $i \in \{1, \cdots ,r\}$, we fix a system $A_i$ (resp. $B_i$) of representatives of right cosets for $Q_i$ (resp. $P_i$) in $S$. Then we define $P$ as the following set of symbols (not as a subset of $G$): $$ P \ := \ S \ \sqcup \ \{ \ at_ia' \ , \ bt_i^{-1}b' \ \vert \ i \in [\![1;r]\!], \ a,b \in S, \ a' \in A_i, \ b' \in B_i \ \} $$

We define $D$ to be the subset of $P \times P$ formed by all the pairs $(x,y)$ listed below. The possible values for parameters in $x$ and $y$ (regarding the above parametrisation) are specified only when some values are not included. We also precise the value of $m(x,y)$ in each case.

\begin{enumerate}[\ding{226}]
\item $(s,s')$, with product $m(s,s') = ss'$ ;
\item $(s,at_ia')$, with product $m(s,at_ia') = (sa)t_ia'$ ;
\item $(s,bt_i^{-1}b')$, with product $m(s,bt_i^{-1}b') = (sb)t_i^{-1}b'$ ;
\item $(at_ia',s)$, with product $m(at_ia',s) = at_i(a's)$, which is rewritten $(au)t_ia''$, for $u \in P_i$ and $a'' \in A_i$ satisfying $\phi_i(u)a'' = a's$ ;
\item $(bt_i^{-1}b',s)$, with product $m(bt_i^{-1}b',s) = bt_i^{-1}(b's)$, which is rewritten $(bv)t_i^{-1}b''$, for $v \in Q_i$ and $b'' \in B_i$ satisfying $\phi_i^{-1}(v)b'' = b's$ ;
\item $(at_ia',bt_i^{-1}b')$ if and only if $a'b \in Q_i$, with product $m(at_ia',bt_i^{-1}b') = a\phi_i^{-1}(a'b)b'$ ;
\item $(bt_i^{-1}b',at_ia')$ if and only if $b'a \in P_i$, with product $m(bt_i^{-1}b',at_ia') = b\phi_i(b'a)a'$.
\end{enumerate}

\noindent The inverse operation $x \mapsto x^{-1}$ is defined on $P$ in the following way:
\begin{enumerate}[\ding{252}]
\item $(s)^{-1} = s^{-1} \in S$ ;
\item $(at_ia')^{-1} = a'^{-1}t_i^{-1}a^{-1}$, which is rewritten $(a'^{-1}v)t_i^{-1}b'$, for $v \in Q_i$ and $b' \in B_i$ satisfying $\phi_i^{-1}(v)b' = a^{-1}$ ;
\item $(bt_i^{-1}b')^{-1} = b'^{-1}t_ib^{-1}$, which is rewritten $(b'^{-1}u)t_ia'$, for $h \in P_i$ and $a' \in A_i$ satisfying $\phi_i(u)a' = b^{-1}$.
\end{enumerate}

In the following, we include some implicit hypotheses in our notations. First, unless specified, any letter $x$ appearing as a subscript in $t_x$, $t_x^{-1}$, $P_x$, $Q_x$, $A_x$ or $B_x$ signifies that $x$ is an integer belonging to $\{1,\cdots ,r\}$. When dealing with elements of $P$, the letter $s$ (or one of its variants such as $s'$ or $s_j$ for $j \in \Nb$) stands for an element of $S$ seen as a subset of $P$. Similarly, denoting an element of $P$ by $xt_iy$ (where $x$ and $y$ are some letters) implicitly means that $x$ and $y$ are elements of $S$ such that $y \in A_i$, and denoting an element of $P$ by $xt_i^{-1}y$ (where $x$ and $y$ are some letters) implicitly means that $x$ and $y$ are elements of $S$ such that $y \in B_i$.

\subsubsection{Inclusion of $P$ in $G$}

The chosen parametrisation allows to embed $P$ in $G$ just by sending the elements $s$, $at_ia'$ and $bt_i^{-1}b'$ of $P$ on the corresponding elements of $G$. To prove that this mapping is an inclusion, the easiest way is to use Britton's Lemma for HNN extensions.

\begin{defi} \label{def:redseq HHN simple}
Let $H := U *_\alpha$ be the HNN extension of a group $U$ relative to an isomorphism $\alpha \colon U_1 \to U_2$, with $U_1,U_2 \leq U$. We denote $t$ the stable letter. A sequence $(u_0,t^{\varepsilon_1},u_1, \cdots , t^{\varepsilon_n},u_n)$, where $n \in \Nb$, each $\varepsilon_j$ is in $\{-1,1\}$ and each $u_j$ belongs to $U$ is said to be \emph{reduced} if there is no consecutive subsequence of the form $(t^{-1},u_j,t)$ with $u_j \in U_1$, or $(t,u_j,t^{-1})$ with $u_j \in U_2$.
\end{defi}

\begin{prop}[{\cite[Chapter IV, Britton's Lemma]{LySch}}]\label{prop Britton HNN}
With the notations of Definition \ref{def:redseq HHN simple}, if the sequence $(u_0,t^{\varepsilon_1},u_1, \cdots , t^{\varepsilon_n},u_n)$ is reduced and $n \geq 1$, then $u_0t^{\varepsilon_1}u_1 \cdots  t^{\varepsilon_n}u_n \neq 1$ in $H$.
\end{prop} 

As a consequence, we can get an analogous result for the group $G$, which can be obtained by a succession of HNN extensions from $S$.

\begin{defi} \label{def:redseq HNN multiple}
A sequence $(s_0,t_{i_1}^{\varepsilon_1},s_1, \cdots , t_{i_n}^{\varepsilon_n},s_n)$, where $n \in \Nb$, each $\varepsilon_j$ is in $\{-1,1\}$ and each $s_j$ belongs to $S$, is said to be \emph{reduced} if it admits no consecutive subsequence of the form $(t_i^{-1},s_j,t_i)$ with $s_j \in P_i$ or $(t_i,s_j,t_i^{-1})$ with $s_j \in Q_i$. 
\end{defi}

\begin{cor} \label{cor Britton LS}
With the notations of Definition \ref{def:redseq HNN multiple}, if the sequence $(s_0,t_{i_1}^{\varepsilon_1},s_1, \cdots , t_{i_n}^{\varepsilon_n},s_n)$ is reduced and $n \geq 1$, then $s_0t_{i_1}^{\varepsilon_1}s_1 \cdots  t_{i_n}^{\varepsilon_n}s_n \neq 1$ in $G$.
\end{cor}

\begin{dem}
First, notice that, in Definition \ref{def:redseq HHN simple}, for any $u \in U$, $(u_0,t^{\varepsilon_1},u_1, \cdots , t^{\varepsilon_n},u_n)$ is reduced if and only if $(u^{-1}u_0,t^{\varepsilon_1},u_1, \cdots , t^{\varepsilon_n},u_n)$ is reduced. In particular, if $u_0t^{\varepsilon_1}u_1 \cdots  t^{\varepsilon_n}u_n = u \in U$, then the sequence $(u_0,t^{\varepsilon_1},u_1, \cdots , t^{\varepsilon_n},u_n)$ is necessarily not reduced by Proposition \ref{prop Britton HNN}. It is this particular formulation of the statement that we use in the following proof.
 
Assume $w:=(s_0,t_{i_1}^{\varepsilon_1},s_1, \cdots , t_{i_n}^{\varepsilon_n},s_n)$ is a sequence such that $s_0t_{i_1}^{\varepsilon_1}s_1 \cdots  t_{i_n}^{\varepsilon_n}s_n = 1$. Let $K^{(0)}$ be the set of all $i_j$ for $j \in \{1,\cdots ,n\}$. If $K^{(0)}$ is empty, then necessarily $n=0$ and we are done. Otherwise, we will prove that $w$ is not reduced. Pick $k_1 \in K^{(0)}$. We can see $G$ as an HNN extension of a certain group $G^{(1)}$ relative to $\phi_{k_1}$, $G^{(1)}$ being the ``HNN extension tower'' of $S$ relative to $\phi_i$ for every $i \in \{1, \cdots ,r\} \setminus \{k_1\}$. Now we can reduce $(s_0,t_{i_1}^{\varepsilon_1},s_1, \cdots , t_{i_n}^{\varepsilon_n},s_n)$ to get a word in $G^{(1)}$ and the symbols $t_{k_1}$, $t_{k_1}^{-1}$, and then apply Proposition \ref{prop Britton HNN} in the HNN extension $G = G^{(1)} *_{\phi_{k_1}}$. It implies that there exists a subsequence $w^{(1)}$ of $w$, composed only with elements of $S$ and symbols $t_k$, $t_k^{-1}$, for $k \in K^{(0)}\setminus \{k_1\}$, whose product belongs either to $P_{k_1}$, in which case $(t_{k_1}^{-1},w^{(1)},t_{k_1})$ is a subsequence of $w$, or to $Q_{k_1}$, in which case $(t_{k_1},w^{(1)},t_{k_1}^{-1})$ is a subsequence of $w$. Let us denote $K^{(1)}$ the subset of $K^{(0)}$ containing the index of the stable letters appearing in $w^{(1)}$. The size of $K^{(1)}$ is strictly less than that of $K^{(0)}$.

We prove by induction on $j \in \Nb^*$ that either $w$ is not reduced, or there exists a tuple $(k_j,K^{(j)},w^{(j)})$ where $k_j \in K^{(j-1)}$, $K^{(j)}$ is a proper subset of $K^{(j-1)}$, and $w^{(j)}$ is a subsequence of $w$ containing only letters in $S$ or symbols $t_k$, $t_k^{-1}$ with $k \in K^{(j)}$, satisfying either that the product of $w^{(j)}$ is in $P_{k_j}$ and $(t_{k_j}^{-1},w^{(j)},t_{k_j})$ is a subsequence of $w$, or that the product of $w^{(j)}$ is in $Q_{k_j}$ and $(t_{k_j},w^{(j)},t_{k_j}^{-1})$ is a subsequence of $w$. 

The case $j=1$ is treated above. Now assume we already proved case $j \in \Nb$, and let us prove case $j+1$. If $K^{(j)}$ is empty, we are done because it means that $w^{(j)}$ is just a letter in $S$, either in $P_{k_j}$ with $(t_{k_j}^{-1},w^{(j)},t_{k_j})$ contained in $w$, or in $Q_{k_j}$ with $(t_{k_j},w^{(j)},t_{k_j}^{-1})$ contained in $w$, and in both cases $w$ is not reduced. 

Otherwise, consider any element $k_{j+1} \in K^{(j)}$. We can see $G$ as an HNN extension of a certain group $G^{(j+1)}$ relative to $\phi_{k_{j+1}}$, $G^{(j+1)}$ being the ``HNN extension tower'' of $S$ relative to $\phi_i$ for every $i \in \{1, \cdots ,r\} \setminus \{k_{j+1}\}$. Then we can reduce $w^{(j)}$ to get a word in $G^{(j+1)}$ and the symbols $t_{k_{j+1}}$, $t_{k_{j+1}}^{-1}$, and then apply Proposition \ref{prop Britton HNN} in the HNN extension $G = G^{(j+1)} *_{\phi_{k_{j+1}}}$ to this word, implying that it is not reduced (because the product of $w^{(j)}$ belongs to $S$). Thus there exists a subsequence $w^{(j+1)}$ of $w^{(j)}$, composed only with elements of $S$ and symbols $t_k$, $t_k^{-1}$, for $k \in K^{(j)} \setminus \{k_{j+1}\}$, whose product belongs either to $P_{k_{j+1}}$, in which case $(t_{k_{j+1}}^{-1},w^{(j+1)},t_{k_{j+1}})$ is a subsequence of $w$, or to $Q_{k_{j+1}}$, in which case $(t_{k_{j+1}},w^{(j+1)},t_{k_{j+1}}^{-1})$ is a subsequence of $w$. Denoting by $K^{(j+1)}$ the (proper) subset of $K^{(j)}$ formed by the index of the stable letters appearing in $w^{(j)}$, this concludes the induction.

The size of $K^{(j)}$ strictly decreases as $j$ grows, but those are finite sets, so the process has to stop and we necessarily get that $w$ is not reduced.
\end{dem}

Now if $s$ and $ct_i^\varepsilon c'$ in $P$ satisfy $s = ct_i^\varepsilon c'$ in $G$, we can rewrite this equality as $s^{-1}ct_i^\varepsilon c' =1$, and Britton's Lemma to get a contradiction (because $(s^{-1}c,t_i^\varepsilon ,c')$ obviously is a reduced word). It is also clear that two elements of $S$ which are distinct in $P$ can't be equal in $G$. Which leaves us with the case where two elements $ct_i^{\varepsilon} c'$ and $dt_j^{e} d'$ in $P$ are equal in $G$. This means that $ct_i^\varepsilon c'd'^{-1}t_j^{-e}d^{-1} = 1$. We can apply Corollary \ref{cor Britton LS} to the word $(c,t_i^\varepsilon ,c'd'^{-1},t_j^{-e},d^{-1})$, which implies that necessarily $j=i$, $e=\varepsilon$ and either $\varepsilon=1$ and $c'd'^{-1} \in Q_i$, or $\varepsilon=-1$ and $c'd'^{-1} \in P_i$. In both cases, as we chose $c'$ and $d'$ to be fixed representatives of right cosets for $Q_i$ or $P_i$ in $S$, we get that $c'=d'$. Now $ct_i^\varepsilon c'd'^{-1}t_j^{-e}d^{-1} = 1$ becomes $cd^{-1} = 1$, so $c=d$ and $ct_i^{\varepsilon} c'$ and $dt_j^{e} d'$ are equal in $P$.

\begin{rem}
Corollary \ref{cor Britton LS} can also be used to prove a normal form theorem for ``HNN extension towers'' similar to $G$.
\end{rem}

Thus, there is a natural inclusion of $P$ in $G$. Moreover, one can check that whenever $(x,y) \in D$, the product $m(x,y)$ in $P$ coincides with the multiplication $xy$ in $G$ through this inclusion, and the inverses of an element $x$ in $P$ and $G$ also coincide.

\subsubsection{Proof that $P$ is a pregroup}

\begin{lem}\label{lem property of P}
For any $x,y \in P$ and $s \in S \subseteq P$, we have $(s,x) \in D$, and $(x,y) \in D$ if and only if $(sx,y) \in D$. Similarly, we have $(y,s) \in D$, and $(x,y) \in D$ if and only if $(x,ys) \in D$.
\end{lem}

\begin{dem}
First, $D$ contains all elements of the form $(s,x)$ for $x \in P$ and $s \in S$. To prove the equivalence $(x,y) \in D \iff (sx,y) \in D$, as multiplying by an element of $S$ does not change the ``type'' of $x$ (element of $S$ ; $at_ia'$ ; or $bt_i^{-1}b'$), we only have to check a few cases. Since the result is clear if $x \in S$ or $y \in S$, we can assume that $x = at_ia'$ and $y = bt_i^{-1}b'$ (or the converse, which is similar). In this case, $$ (x,y) \in D \iff a'b \in Q_i \iff \big((sa)t_ia',bt_i^{-1}b' \big) \in D \iff (sx,y) \in D \ . $$
The proof of the other assertion is similar.
\end{dem}

\begin{prop}
As defined above, $(P,D)$ is a pregroup.
\end{prop}

\begin{dem}
Axioms \ref{cond:Pre1} and \ref{cond:Pre2} are easily verified. For Axiom \ref{cond:Pre3}, we only have to check the conditions on $D$, because the associativity of the product is a consequence of the fact that multiplication in $P$ and $G$ coincide. Considering Lemma \ref{lem property of P}, we are left with only two main cases to check for Axiom \ref{cond:Pre3}. Here, $X_i$ stands for $t_i$ or $t_i^{-1}$ (and $X_i^{-1}$ stands for the other one). 
\begin{enumerate}[-]
\item If $(cX_ic',s) \in D$ and $(s,dX_jd') \in D$, then we have $(cX_ic's,dX_jd') \in D \iff \big(cX_ic',(sd)X_jd' \big) \in D$ because the condition ($i=j$, $(X_i)^{-1} = X_j$ and $c'sd \in P_i$ or $Q_i$) is the same in both cases.
\item If $(cX_ic',dX_i^{-1}d') \in D$ and $(dX_i^{-1}d',eX_ie') \in D$: the product of each of these pairs belongs to $S$, so everything is defined.
\end{enumerate}

Finally, for Axiom \ref{cond:Pre4}, take $(w,x),(x,y),(y,z) \in D$. If $w \in S$, $z \in S$ or $m(x,y) \in S$, then the conclusion holds. Otherwise, we necessarily have $x \in S$ or $y \in S$, so $(x,m(y,z)) \in D$ or $(m(w,x),y) \in D$, and Axiom \ref{cond:Pre3} allows to conclude.
\end{dem}

\begin{cor}
The finite pregroup $P$ has universal group $G$, the Leary-Stancu group associated to $\F$ and the family of generators $\{ \phi_1, \cdots , \phi_r \}$.
\end{cor}

\begin{dem}
We already proved that $P$ is a subpregroup of $G$, because it is a pregroup contained in $G$ and the multiplication laws are compatible. Moreover, $P$ contains $S$ and the $t_i$ for $i \in \{1,\cdots ,r\}$, which generate $G$, so $G = U(P)$ by Proposition \ref{prop:universalgroup}.
\end{dem}

\subsubsection{Can $P$ be a locality?}

Here we ask whether the finite pregroup $P$ constructed above (or more precisely: the underlying partial group of $P$) can be a locality over $S$, for a certain set of objets $\Delta$.

First, let us describe the domain $\Db_P$ of $P$ when it is considered as a partial group.

\begin{prop}
The domain $\Db_P$ is constituted of all the words $w \in \Wb(P)$ satisfying the following conditions:
\begin{enumerate}[$(i)$]
\item\label{cond:domain1} $w$ does not contain simultaneaously a term of the form $ct_i^\varepsilon c'$ and another of the form $dt_j^{\varepsilon'}d'$ unless $i = j$. Hence $w$ only contains terms of the form $s \in S$, $at_ia'$ or $bt_i^{-1}b'$ for a fixed $i$.
\item\label{cond:domain2} The terms in $w$ that are not elements of $S$ should alternate between the forms $at_ia'$ and $bt_i^{-1}b'$ (possibly with terms in $S$ interposed).
\item\label{cond:domain3} Between any term $at_ia'$ in $w$ and the next term of the form $bt_i^{-1}b'$, if we denote by $s$ the product of all the (possible) terms in $S$ interposed between $at_ia'$ and $bt_i^{-1}b'$, then we should have $a'sb \in Q_i$.
\item\label{cond:domain4} Between any term $bt_i^{-1}b'$ in $w$ and the next term of the form $at_ia'$, if we denote by $s$ the product of all the (possible) terms in $S$ interposed between $bt_i^{-1}b'$ and $at_ia'$, then we should have $b'sa \in Q_i$.
\end{enumerate}
\end{prop}

\begin{dem}[sketch]
Recall that we defined the domain $\Db_P$ as follows (just before Proposition \ref{prop pregroup is partial group}): $$\Db_P \ = \ \left\{ \ (x_1,x_2,\dots,x_n) \in \Wb(P) \ \vert \ \forall k,l\in\{1,2,\dots, n\} \textup{ with } k<l, \ \ x_kx_{k+1}\cdots x_l\in P \ \right\} $$

First, let $w = (x_1,x_2,\dots,x_n) \in \Wb(P)$ be any word satisfying the four conditions in the statement. Denote by $i$ the index of the symbol $t_i$ possibly appearing in the terms of $w$ (unique by condition \ref{cond:domain1}). We can prove by induction on $m \in \Nb$ that for any subword $(x_k,\cdots , x_{k+m})$ of $w$, we are in one of the following cases (where the symbols $y_1, \cdots,y_q$ refer to any allowed elements of $P$): \begin{enumerate}[-]
\item $x_kx_{k+1}\cdots x_{k+m} = at_ia'$ if $(x_k,\cdots , x_{k+m})= (s_1,\cdots , s_j,ct_ic',y_1,\cdots ,y_q , dt_id',s_1',\cdots , s_{j'}')$ where $s_1,\cdots , s_j,s_1',\cdots , s_{j'}' \in S$ (with $j,j',q \in \Nb$), and in this case we have $a \in (s_1\cdots  s_jc)P_i$ and $a' \in Q_i(d's_1' \cdots s_{j'}')$ ;
\item $x_kx_{k+1}\cdots x_{k+m} = bt_i^{-1}b'$ if $(x_k,\cdots , x_{k+m})= (s_1,\cdots , s_j,ct_i^{-1}c',y_1,\cdots ,y_q , dt_i^{-1}d',s_1',\cdots , s_{j'}')$ \linebreak[4] where $s_1,\cdots , s_j,s_1',\cdots , s_{j'}' \in S$ (with $j,j',q \in \Nb$), and in this case we have $b \in (s_1\cdots  s_jc)Q_i$ and $b' \in P_i(d's_1' \cdots s_{j'}')$ ;
\item $x_kx_{k+1}\cdots x_{k+m}=s \in S$ if $(x_k,\cdots , x_{k+m})= (s_1,\cdots , s_j,ct_i^{-1}c',y_1,\cdots ,y_q , dt_id',s_1',\cdots , s_{j'}')$ where $s_1,\cdots , s_j,s_1',\cdots , s_{j'}' \in S$ (with $j,j',q \in \Nb$), and in this case $s \in (s_1\cdots  s_jc)Q_i(d's_1' \cdots s_{j'}')$ ;
\item $x_kx_{k+1}\cdots x_{k+m}=s \in S$ if $(x_k,\cdots , x_{k+m})= (s_1,\cdots , s_j,ct_ic',y_1,\cdots ,y_q , dt_i^{-1}d',s_1',\cdots , s_{j'}')$ where $s_1,\cdots , s_j,s_1',\cdots , s_{j'}' \in S$ (with $j,j',q \in \Nb$), and in this case $s \in (s_1\cdots  s_jc)P_i(d's_1' \cdots s_{j'}')$ ;
\item $x_kx_{k+1}\cdots x_{k+m}=s \in S$ if $(x_k,\cdots , x_{k+m})= (s_1,\cdots , s_{m+1})$ where $s_1,\cdots , s_{m+1} \in S$, and in this case $s = s_1 \cdots s_{m+1}$.
\end{enumerate}
The case $m=0$ is trivial. Assume the fact is proved for the rank $m \in \Nb$. In order to prove it for the rank $m+1$, take a subword $(x_k,\cdots , x_{k+m+1})$ of $w$, and apply the induction hypothesis on $(x_k,\cdots , x_{k+m})$. Then check what we need for $(x_k,\cdots , x_{k+m+1})$, distinguishing between the five cases listed above for $(x_k,\cdots , x_{k+m})$, also distinguishing between the two or three possibilities for the term $x_{k+m+1}$, and using conditions \ref{cond:domain3} and \ref{cond:domain4}.

Now let us prove by contraposition that any word in $\Db_P$ satisfies the four conditions. If $w \in \Wb(P)$ does not satisfy condition \ref{cond:domain1}, consider the smallest subword of $w$ of the form $(ct_i^\varepsilon c', x_1, \cdots ,x_k,dt_j^{\varepsilon'}d')$ with $i \neq j$. Necessarily, $x_1, \cdots ,x_k \in S$ so $x_1 \cdots x_k =: s \in S$, so the product of $(ct_i^\varepsilon c', x_1, \cdots ,x_k,dt_j^{\varepsilon'}d')$ belongs to $P$ if and only if $(ct_i^\varepsilon c's,dt_j^{\varepsilon'}d') \in D$, which is excluded since $i \neq j$.

If $w \in \Wb(P)$ does not satisfy condition \ref{cond:domain2}, we can consider the smallest subword of $w$ of the form $(ct_i^\varepsilon c', x_1, \cdots ,x_k,dt_i^\varepsilon d')$. Necessarily, $x_1, \cdots ,x_k \in S$ so $x_1 \cdots x_k =: s \in S$, so the product of the sequence $(ct_i^\varepsilon c', x_1, \cdots ,x_k,dt_i^\varepsilon d')$ belongs to $P$ if and only if $(ct_i^\varepsilon c's,dt_i^{\varepsilon}d') \in D$, which can't be true.

If $w \in \Wb(P)$ does not satisfy condition \ref{cond:domain3}, then we can consider a subword of $w$ of the form $(at_ia', x_1, \cdots ,x_k,bt_i^{-1}b')$ where each $x_j$ belongs to $S$ and $a'x_1 \cdots x_kb \notin Q_i$. Then $\big(at_i(a'x_1 \cdots x_k),bt_i^{-1}b' \big)$ is not in $D$, so the product $(at_ia', x_1, \cdots ,x_k,bt_i^{-1}b')$ is not in $P$. The argument is identical if $w \in \Wb(P)$ does not satisfy condition \ref{cond:domain4}.
\end{dem}

Now we can look for a candidate for the set of objects $\Delta$, included in the set of subgroups of $S$ since we want $(P,\Delta,S)$ to be a locality. Remark that for any subgroup $R \leq S$, we have $R \subseteq \Db_P(s)$ for all $s \in S$ (using Notation \ref{notation conjugaison}). Moreover, for all $i \in \{1, \cdots ,r\}$, $a \in S$ and $a' \in A_i$ we have \begin{align*}
R \subseteq \Db_P(at_ia') &\iff \forall s \in R, \quad \big((at_ia')^{-1}s,at_ia' \big) \in D \\
&\iff \forall s \in R, \quad a^{-1}sa \in P_i \\
&\iff R \leq {}^a P_i
\end{align*}
In this case, we get $R^{at_ia'} = {\big(\phi_i(R^a)\big)}^{a'}$. Similarly, $R \subseteq \Db_P(bt_i^{-1}b') \iff R \leq {}^b Q_i$ and in that case $R^{bt_i^{-1}b'} = {\big(\phi_i^{-1}(R^b)\big)}^{b'}$.

Since $\Db_P$ contains all the words of length one $(at_ia')$ and $(bt_i^{-1}b')$, since $\Delta$ has to be closed under taking overgroups in $S$, and because we need $\Db_P = \Db_\Delta$, the above remark implies that $\Delta$ must contain the groups $P_i$ and $Q_i$ and all their $S$-conjugate for every $i \in \{1,\cdots ,r\}$.

However, the equality $\Db_P = \Db_\Delta$ can't hold whenever there exists $i \neq j$ such that $P_i = P_j$. Indeed, in such a case $(t_i^{-1},t_j)$ is in $\Db_\Delta$ via $(Q_i,P_i,Q_j)$, but $(t_i^{-1},t_j) \notin \Db_P$. Similarly, we can't have $\Db_P = \Db_\Delta$ if $P_i = Q_i$ for a certain $i$, because in this case $(t_i,t_i)$ is in $\Db_\Delta$ via $(P_i,P_i,P_i)$ but not in $\Db_P$.

Even worse, if there is an $i \in \{1, \cdots ,r\}$ such that $P_i < S$, then there exists $s \in N_S(P_i) \setminus P_i$ and $(t_i^{-1},s,t_i)$ is in $\Db_\Delta$ via $(Q_i,P_i,P_i,Q_i)$, but not in $\Db_P$. In conclusion, the pregroup $P$ is never a locality.

\subsection{A pregroup for the Robinson group}

Let $\F$ be a fusion system over $S$ generated by a family $\{\F_{S_1}(G_1), \cdots , \F_{S_r}(G_r) \}$ of realisable fusion subsystems, where $S_1, \cdots , S_r$ are subgroups of $S$ and $G_1,\cdots ,G_r$ are finite groups. For each $i \in [\![1;r]\!]$, assume that $S_i$ is contained in $G_i$ as a Sylow $p$-subgroup via a morphism $f_i \colon S_i \hookrightarrow G_i$, whose image will be denoted by $S_i'$. Let $F$ be the free product of $S$ and the groups $G_i$. Then the Robison group $G$ is defined to be the quotient of $F$ by the normal closure of the elements $uf_i(u)^{-1}$ for any $i \in \{1,\cdots ,r\}$ and $u \in S_i$. The group $G$ can be seen as an iterated free product with amalgamation: $$G = \Big( \cdots \big( (S *_{S_1} G_1) *_{S_2} G_2 \big) \cdots \Big) *_{S_r} G_r $$

\begin{rem}
In his article \cite{Ro}, Robinson states a result (Theorem 2) saying that $\F = \F_S(G)$, but only in the case of ``Alperin fusion systems'', which is a class of fusion systems containing saturated fusion systems. However, Theorem 1 in the same article implies that the result remains true for the larger class of fusion systems generated by families of realisable fusion subsystems, which is our framework here.
\end{rem}

\subsubsection{Constructing the pregroup}

Informally, $P$ is the subset of $G$ containing $S$ together with all the elements of the form $aga'$ where $g \in G_i$ for a certain $i \in \{1, \cdots , r\}$ and $a,a' \in S$. A pair $(x,y)$ belongs to the domain $D$ if and only if the product $xy$ in $G$ belongs to $P$, and multiplication and inverses of elements of $P$ are then defined as in $G$.

Because of the relations that exist in $G$, elements of the form $aga'$ in $G$ can admit several such representations. In order to get a bijective parametrisation, we need to fix for each $i \in \{1, \cdots ,r\}$ a system $H_i$ of representatives of \emph{non-trivial} double cosets $S_i'gS_i'$ in $G_i$. Moreover we also need to introduce, for each $g \in H_i$, the subgroup $T_{i,g}:= f_i^{-1} \big(S_i' \cap {}^gS_i' \big)$ of $S_i$, and fix a system $A_{i,g}$ of representatives of cosets in $S/T_{i,g}$. Beware that $T_{i,g}$ depends on the choice of the representative in the double coset $S_i'gS_i'$ (however, the $S_i$-conjugacy class of $T_{i,g}$ does not depend on it). Then we can define $P$ explicitly as the following set of symbols: $$ P \ := \ S \ \sqcup \ \{ \ aga' \ \vert \ i \in [\![1;r]\!], \ g \in H_i, \ a \in A_{i,g}, \ a' \in S \ \} $$

We define $D$ to be the subset of $P \times P$ formed by all the pairs $(x,y)$ listed below. The possible values for parameters in $x$ and $y$ (regarding the above parametrisation) are specified only when some values are not included. We also precise the value of $m(x,y)$ in each case.
\begin{enumerate}[\ding{226}]
\item $(s,s')$, with product $m(s,s') = ss'$ ;
\item $(s,aga')$, with product $m(s,aga') = (sa)ga'$, which is rewritten $a''g \big( f_i^{-1}(g^{-1}f_i(v)g)a' \big)$ for $a'' \in A_{i,g}$ and $v \in T_{i,g}$ satisfying $a''v = sa$ (where $i \in \{1,\cdots,r\}$ is such that $g \in H_i$) ;
\item $(aga',s)$, with product $m(aga',s) = ag(a's)$ ;
\item $(aga',bhb')$ if and only if $g,h \in H_i$, for a certain $i \in \{1, \cdots ,r\}$, and $a'b \in S_i$. The product $m(aga',bhb')$ then equals $a {f_i}^{-1}\big(gf_i(a'b)h\big)b' \in S$ if $gf_i(a'b)h \in G_i$ belongs to $S_i'$. Otherwise, we have $gf_i(a'b)h = f_i(u)kf_i(u')$ for certain $k \in H_i$ and $u,u' \in S_i$, and the product $m(aga',bhb')$ equals $(au)k(u'b')$, which is rewritten $ck\big( f_i^{-1}(k^{-1}f_i(v)k)u'b' \big)$ for $c \in A_{i,k}$ and $v \in T_{i,k}$ satisfying $cv = au$.
\end{enumerate}

\noindent The inverse operation $x \mapsto x^{-1}$ is defined on $P$ in the following way:
\begin{enumerate}[\ding{252}]
\item $(s)^{-1} = s^{-1} \in S$ ;
\item $(aga')^{-1} = a'^{-1}g^{-1}a^{-1}$ in $G$, which equals $(a'^{-1}u)h(u'a^{-1})$ for certain $h \in H_i$ (index $i$ being the same as for $g$) and $u,u' \in S_i$ (satisfying $f_i(u)hf_i(u') = g^{-1}$ in $G_i$), and leads to the element $a''h \big( f_i^{-1}(h^{-1}f_i(v)h)u'a^{-1} \big)$ in $P$, with $a'' \in A_{i,h}$ and $v \in T_{i,h}$ satisfying $a''v = a'^{-1}u$.
\end{enumerate}

In the following, we include some implicit hypotheses in our notations. First, unless specified, any letter $x$ appearing as a subscript in $f_x$, $G_x$, $S_x$, $S_x'$ or $H_x$ signifies that $x$ is an integer belonging to $\{1,\cdots ,r\}$. When dealing with elements of $P$, the letter $s$ (or one of its variants such as $s'$ or $s_j$ for $j \in \Nb$) stands for an element of $S$ seen as a subset of $P$. Similarly, denoting an element of $P$ by $xgy$ when $g$ is in some $H_i$ (and with $x$ and $y$ being some letters) implicitly means that $x$ and $y$ are elements of $S$ such that $x \in A_{i,g}$.

\begin{rem}
In the case where $r=1$, we already know a finite pregroup whose universal group is $S *_{S_1} G_1$, namely $S \cup_{S_1} G_1$ as in Example \ref{ex:A*_CB}. Notice that this pregroup does not coincide with the set $P$ we just introduced. In fact, it is strictly contained in $P$. One could ask if instead of defining $P$ as above, we could just define it to be $\Big( \cdots (S \cup_{S_1} G_1) \cup_{S_2} G_2 \cdots \Big) \cup_{S_r} G_r$. This is not a pregroup in general as soon as $r \leq 2$, because of the need for Axiom \ref{cond:Pre4} to be satisfied.
\end{rem}

\subsubsection{Inclusion of $P$ in $G$}

We defined the elements of $P$ as symbols of certain elements in $G$, so there is an obvious mapping $P \to G$. Moreover, the products of elements in $P$ (when defined) coincide with the products of the same elements in $G$. Indeed, the only cases where it is not obvious in the way we defined $m(x,y)$ is when it is rewritten, and all the rewritings result from the relations in $G$ such as $f_i(u) =u$ for all $u \in S_i$, or $vg = g f_i^{-1}(g^{-1}f_i(v)g)$ for all $g \in H_i$ and $v \in T_{i,g}$.

There remains to prove that two distinct elements of $P$ can not be equal when considered as elements of $G$. For this, we need the Normal Form Theorem for free product with amalgamation, which we state below.

\begin{defi} \label{def:redseq amal}
Let $A *_C B$ be the free product of groups $A$ and $B$ amalgamating $C$, where $C$ is a subgroup of $A$ isomorphic to a subgroup $C'$ of $B$ via a morphism $f$. A sequence $(x_1,\cdots ,x_n)$ of elements of $A *_C B$ is said to be \emph{reduced} if:
\begin{enumerate}[-]
\item each $x_i$ belongs to one of the factors $A$ or $B$ ;
\item $x_i$ and $x_{i+1}$ always come from different factors ;
\item if $n>1$, no $x_i$ belongs to $C$ or $C'$ ;
\item if $n=1$, then $x_1 \neq 1$.
\end{enumerate}
\end{defi}

\begin{prop}[{\cite[Chapter IV, Theorem 2.6]{LySch}}] \label{prop NFT FPA}
With the notations of Definition \ref{def:redseq amal}, if $(x_1,\cdots ,x_n)$ is reduced, then $x_1\cdots x_n \neq 1$ in $A *_C B$.
\end{prop}

First, if $s$ and $aga'$ in $P$, with $g$ in some $H_i$, are equal as elements of $G$, then $s^{-1}aga'=1$ in $S *_{S_i} G_i$ (seen as a subgroup of $G$). Applying Proposition \ref{prop NFT FPA} (more precisely its contrapositive) to the sequence $(s^{-1}a,g,a')$, we get that $s^{-1}a \in S_i$ or $a' \in S_i$. Replacing $g$ with $f_i(s^{-1}a)g$, $gf_i(a')$ or $f_i(s^{-1}a)gf_i(a')$, we can apply Proposition \ref{prop NFT FPA} to a new reduced sequence whose product is $s^{-1}aga'$ and get a contradiction, so $s$ and $aga'$ can't be equal in $G$.

Now assume that $aga'$ and $bhb'$ are two elements of $P$, with $g$ in some $H_i$ and $h$ in some $H_j$, such that $aga' = bhb'$ in $G$, i.e. $a'^{-1}g^{-1}a^{-1}bhb' = 1$. If $i \neq j$, we can consider this equality in $G_i *_{S_i} \big(S *_{S_j} G_j \big)$ (seen as a subgroup of $G$) and apply Proposition \ref{prop NFT FPA} to the sequence $(a'^{-1},g^{-1},a^{-1}bhb')$, leading to a contradiction, unless $a'^{-1} \in S_i$ or $a^{-1}bhb' \in S_i$ (or both). Up to replacing $g^{-1}$ with another representative of its double coset, we can apply Proposition \ref{prop NFT FPA} to a reduced sequence and get a contradiction.

Hence $i=j$, and we want to apply Proposition \ref{prop NFT FPA} to $(a'^{-1},g^{-1},a^{-1}b,h,b')$ seen as a sequence in $S *_{S_i} G_i$. If $a'^{-1}$ or $b'$ is in $S_i$, we can do the same trick as before to get a reduced sequence. Proposition \ref{prop NFT FPA} then implies that $a^{-1}b$ is in $S_i$. Thus we can rewrite $a'^{-1}g^{-1}a^{-1}bhb' = a'^{-1} \big(g^{-1}f_i(a^{-1}b)h \big) b'$, and again by Proposition \ref{prop NFT FPA} we get that $g^{-1}f_i(a^{-1}b)h \in S_i'$, and still $a'^{-1} f_i^{-1}\big(g^{-1}f_i(a^{-1}b)h \big) b' = 1$ in $S$ (and these two facts precisely hold with our former $g$ and $h$, no matter the ``tricks'' we had to do). In particular, there exists $u \in S_i'$ such that $gu = f_i(a^{-1}b)h$. As $g$ and $h$ were fixed representatives of double cosets for $S_i'$, this means that $g=h$ and, denoting $v := f_i(a^{-1}b) \in S_i'$, that $g^{-1}vg \in S_i'$. This amounts to $v \in S_i' \cap {}^g S_i'$, or equivalently $a^{-1}b \in T_{i,g}$. Thus $b \in aT_{i,g}$, but we chose $a$ and $b$ to be representatives of left cosets for $T_{i,g}$ in $S_i$, so $a =b$. Now $a'^{-1} f_i^{-1}\big(g^{-1}f_i(a^{-1}b)h \big) b' = 1$ becomes $a'^{-1}b'=1$, i.e. $a'=b'$. Finally, $aga'$ and $bhb'$ are equal in $P$, so the natural mapping $P \to G$ is an inclusion.

\subsubsection{Proof that $P$ is a pregroup}

\begin{lem} \label{lem properties of P}
The $(P,D)$ constructed above satisfies the two following properties:
\begin{enumerate}[$(1)$]
\item For all $x,y \in P$ and $s \in S \subseteq P$, we always have $(s,x) \in D$, and $(x,y) \in D$ if and only if $(sx,y) \in D$. Similarly, we always have $(y,s) \in D$, and $(x,y) \in D$ if and only if $(x,ys) \in D$.
\item Let $y \in P \setminus S$. If $x \in P$ satisfies $(x,y) \in D$, then for all $z \in P$, $(y,z) \in D$ implies $(xy,z) \in D$. Similarly, if $z \in P$ satisfies $(y,z) \in D$, then for all $x \in P$, $(x,y) \in D$ implies $(x,yz) \in D$.
\end{enumerate}
\end{lem}

\begin{dem}
For the first property, let $x,y \in P$ and $s \in S \subseteq P$. The fact that $(s,x) \in D$ (and $(y,s) \in D$) is clear in how we defined $D$. Now if $y \in S$, the equivalence $(x,y) \in D \iff (sx,y) \in D$ is true for the same reason, and it is also true if $x \in S$ (because then $sx \in S$ too). So we can assume $x = aga'$ and $y=bhb'$ with $g$ in some $H_i$, $h$ in some $H_j$, $a \in A_{i,g}$, $b \in A_{j,h}$ and $a',b' \in S$. Then $sx = a''g \big( f_i^{-1}(g^{-1}f_i(v)g)a' \big)$ for $a'' \in A_{i,g}$ and $v \in T_{i,g}$ satisfying $a''v = sa$. Thus \begin{align*}
(x,y) \in D \iff (aga',bhb') \in D &\iff i=j \quad \textup{and} \quad a'b \in S_i \\
&\iff i=j \quad \textup{and} \quad f_i^{-1}(g^{-1}f_i(v)g)a'b \in S_i \\
&\iff (sx,y) \in D \ .
\end{align*}
The proof of $(x,y) \in D \iff (x,ys) \in D$ is similar.

For the second property, fix $y \in P \setminus S$ and $x \in P$ such that $(x,y) \in D$, and take any $z \in P$. If $x \in S$, we have an equivalence $(y,z) \in D \iff (xy,z) \in D$, which comes from the first property. If $z \in S$, we also have an obvious equivalence. So we can assume that $x = aga'$ with $g$ in some $H_i$, $y = bhb'$ with $h$ in the same $H_i$ (because $(x,y) \in D$), and $z = ckc'$ with $k$ in some $H_j$. Moreover we have $a'b \in S_i$ and we can assume that $xy \in P \setminus S$ (otherwise it is clear that $(xy,z) \in D$), i.e. $gf_i(a'b)h = f_i(u)\widetilde{g}f_i(u')$ for certain $\widetilde{g} \in H_i$ and $u,u' \in S_i$. Then $xy = d\widetilde{g}\big( f_i^{-1}(\widetilde{g}^{-1}f_i(v)\widetilde{g})u'b' \big)$ for $d \in A_{i,\widetilde{g}}$ and $v \in T_{i,\widetilde{g}}$ satisfying $dv = au$. Thus \begin{align*}
(y,z) \in D \iff (bhb',ckc') \in D &\iff i=j \quad \textup{and} \quad b'c \in S_i \\
&\iff i=j \quad \textup{and} \quad f_i^{-1}(\widetilde{g}^{-1}f_i(v)\widetilde{g})u'b'c \in S_i \\
&\iff \Big(d\widetilde{g}\big( f_i^{-1}(\widetilde{g}^{-1}f_i(v)\widetilde{g})u'b' \big) \ , \ ckc' \Big) \in D \\
&\iff (xy,z) \in D \ .
\end{align*}
In particular, the implication $(y,z) \in D \implies (xy,z) \in D$ holds. The proof of the last assertion is similar.
\end{dem}

\begin{prop}
As defined above, $(P,D)$ is a pregroup.
\end{prop}

\begin{dem}
Axioms \ref{cond:Pre1} and \ref{cond:Pre2} are easily verified. In order to verify the ``domain part'' of Axiom \ref{cond:Pre3} on a triplet $(x,y,z)$ of elements of $P$, the first property in Lemma \ref{lem properties of P} implies that we can assume $x \notin S$, because otherwise we would have $(xy,z) \in D$ and also $(x,yz) \in D$ (since $x$ would be in $S$). Similarly, we can assume $z \notin S$. By the second property, if $y \notin S$ then we have both $(xy,z)$ and $(x,yz)$ in $D$. Thus we can assume $y \in S$, so that $(x,y,z) = (aga',s,bhb')$ with $g$ in some $H_i$ and $h$ in some $H_j$. In this case we have \begin{align*}
 (xy,z) \in D &\iff \big((aga')s,bhb' \big) \in D \\
 &\iff i=j \quad \textup{and} \quad a'sb \in S_i \\
 &\iff i=j \quad \textup{and} \quad \forall v \in T_{i,h}, \quad a'sbv^{-1} \in S_i \\
 &\iff \big(aga',s(bhb') \big) \in D \qquad \textup{because } s(bhb') = b''h \big( f_j^{-1}(g^{-1}f_j(v)h)b' \big) \textup{ for some} \\ 
 &\phantom{\iff \big(aga',s(bhb') \big) \in D\qquad \textup{because }} \textup{$b'' \in A_{j,h}$ and $v \in T_{j,h}$ satisfying $b'' = sbv^{-1}$.} \\
 &\iff (x,yz) \in D \ .
\end{align*}

One also have to check that the $(xy,z)$ and $(x,yz)$ coincide, but this is a consequence of the fact that $P$ is contained in $G$ with consistent products and the associativity of the group law in $G$.

Finally, for Axiom \ref{cond:Pre4}, assume that $(w,x),(x,y),(y,z) \in D$. The conclusion is then obvious if $w \in S$ or $z \in S$. By the first property in Lemma \ref{lem properties of P}, the conclusion is also true whenever $x \in S$ or $y \in S$. Thus we are left with checking the result in the case where none of the four elements is in $S$, but then it directly follows from the second property in Lemma \ref{lem properties of P}.
\end{dem}

\begin{rem}
One can check that each $G_i$ embeds in $P$ as a subgroup, via the following mapping: if $x \in G_i$ belongs to $S_i$, send it directly to $x \in S$ in $P$ ; otherwise, write it $x = f_i(a)gf_i(a')$ with $g \in H_i$, $a \in A_{i,g}$ and $a' \in S_i$, and send it to $aga'$ in $P$. With this embedding of $G_i$, $S_i$ and $S_i'$ become indentified in $P$ and $S \cap G_i = S_i$.
\end{rem}

\begin{cor}
The finite pregroup $P$ has universal group $G$, the Robinson group associated to $\F$ and the generating family $\{\F_{S_1}(G_1), \cdots , \F_{S_r}(G_r) \}$.
\end{cor}

\begin{dem}
We already proved that $P$ is a subpregroup of $G$, because it is a pregroup contained in $G$ and the multiplication laws are compatible. Moreover, $P$ contains $S$ and the $G_i$ for $i \in \{1,\cdots ,r\}$, which generate $G$, so $G = U(P)$ by Proposition \ref{prop:universalgroup}.
\end{dem}

\subsubsection{Can $P$ be a locality?}

As for the pregroup associated with the Leary-Stancu group, we can wonder if (the underlying partial group of) $P$ can be equipped with a set of objects $\Delta$ such that $(P,\Delta,S)$ is a locality. First, we describe the domain $\Db_P$ of $P$ seen as a partial group.

\begin{prop}
The domain $\Db_P$ is constituted of all the words $w \in \Wb(P)$ satisfying the following conditions:
\begin{enumerate}[$(i)$]
\item $w$ does not contain simultaneously a term of the form $aga'$ and a term of the form $bhb'$ if $g \in H_i$ and $h \in H_j$ with $i \neq j$.
\item Between any term $aga'$ in $w$ and the next term of the form $bhb'$ (necessarily with $g,h \in H_i$ for a fixed $i$), if we denote by $s$ the product of all the (possible) terms in $S$ interposed between $aga'$ and $bhb'$, then we should have $a'sb \in S_i$.
\end{enumerate}
\end{prop}

\begin{dem}[sketch]
The domain $\Db_P$ is defined to be $$\Db_P \ = \ \left\{ \ (x_1,x_2,\dots,x_n) \in \Wb(P) \ \vert \ \forall k,l\in\{1,2,\dots, n\} \textup{ with } k<l, \ \ x_kx_{k+1}\cdots x_l\in P \ \right\}  \ .$$

Let $w = (x_1,x_2,\dots,x_n) \in \Wb(P)$ be a word satisfying the two conditions in the statement. By the first condition, there exists $i \in \{1,\cdots ,r\}$ such that each $x_j$ is either an element of $S$ or an element of the form $aga'$ with $g \in H_i$. We can prove by induction on $m \in \Nb$ that for any subword $(x_k,\cdots , x_{k+m})$ of $w$, we are in one of the following cases (where the symbols $y_1, \cdots,y_q$ refer to any allowed elements of $P$): \begin{enumerate}[-]
\item $x_kx_{k+1}\cdots x_{k+m}$ equals some $aga'$ if $(x_k,\cdots , x_{k+m})= (s_1,\cdots , s_j,bhb',y_1,\cdots ,y_q , c\widetilde{h} c',s_1',\cdots , s_{j'}')$ where $s_1,\cdots , s_j,s_1',\cdots , s_{j'}' \in S$ (with $j,j',q \in \Nb$), and in this case we have $a \in (s_1\cdots  s_jb)S_i$ and $a' \in S_i(c's_1' \cdots s_{j'}')$ ;
\item $x_kx_{k+1}\cdots x_{k+m}$ equals some $s \in S$ if $(x_k,\cdots , x_{k+m})= (s_1,\cdots , s_j,aga',y_1,\cdots ,y_q , bhb',s_1',\cdots , s_{j'}')$ where $s_1,\cdots , s_j,s_1',\cdots , s_{j'}' \in S$ (with $j,j',q \in \Nb$), and in this case $s \in (s_1\cdots  s_ja)S_i(b's_1' \cdots s_{j'}')$ ;
\item $x_kx_{k+1}\cdots x_{k+m}$ equals some $s \in S$ if $(x_k,\cdots , x_{k+m})= (s_1,\cdots , s_{m+1})$ where $s_1,\cdots , s_{m+1} \in S$, and in this case $s = s_1 \cdots s_{m+1}$.
\end{enumerate}
The case $m=0$ is trivial. Assume the fact is proved for the rank $m \in \Nb$. In order to prove it for the rank $m+1$, take a subword $(x_k,\cdots , x_{k+m+1})$ of $w$, and apply the induction hypothesis on $(x_k,\cdots , x_{k+m})$. Then check what we need for $(x_k,\cdots , x_{k+m+1})$, distinguishing between the three cases listed above for $(x_k,\cdots , x_{k+m})$, also distinguishing between the two possibilities for the term $x_{k+m+1}$, and using the two conditions of the statement.

Reciprocally, take $w \in \Db_P$. Consider any subword of $w$ of the form $(aga',s_1,\cdots,s_k,bhb')$, with $s_1, \cdots , s_k \in S$. Then $s:=s_1 \cdots s_k$ is in $S$ too, and we have $(aga')s_1\cdots s_k(bhb') \in P$ if and only if $\big(ag(a's),bhb' \big) \in D$, because all any product involving an element of $S$ is defined in $P$. In particular, $g$ and $h$ must belong to the same $H_i$, and $a'sb$ must belong to $S_i$. Considering that it has to be true for every subword of $w$ of this form, this proves that $w$ fulfills the two conditions.
\end{dem}

Now we can look for a subset $\Delta$ of the subgroups of $S$ such that $(P,\Delta,S)$ would be a locality. Remark that for any subgroup $R \leq S$, we have $R \subseteq \Db_P(s)$ for all $s \in S$. Moreover, for all $i \in \{1, \cdots ,r\}$, $g \in H_i$, $a \in A_{i,g}$ and $a' \in S$ we have \begin{align*}
R \subseteq \Db_P(aga') &\iff \forall s \in R, \quad \big((aga')^{-1}s,aga' \big) \in D \\
&\iff a^{-1}sa \in S_i \\
&\iff R \leq {}^a S_i \ .
\end{align*}
In this case, $R^{aga'}$ is a subgroup of ${S_i^g}^{a'}$. Moreover, $R^{aga'} \leq S \iff R^{ag} \leq S$, and $R^{ag}$ is a subgroup of $S_i^g$ (hence of $G_i$). But for a subgroup of $G_i$ it is equivalent to be contained in $S$ and to be contained in $S_i$, so if we want $R^{aga'} \leq S$, i.e. $R^{ag} \leq S$, it is equivalent to ask $R^{ag} \leq S_i$, which is rewritten $R^a \leq {}^gS_i$. Since we are in the case where $R \leq {}^a S_i$, this amounts to $R^a \leq S_i \cap {}^gS_i$, i.e. $R \leq {}^a T_{i,g}$.

Thus, because $\Db_P$ contains all the words of length one $(aga')$ and $\Delta$ has to be closed under taking overgroups in $S$, $\Delta$ must contain all the subgroups $T_{i,g}$, for $i \in \{1, \cdots ,r\}$ and $g \in H_i$, together with all their $S$-conjugate (and the overgroups of all these in $S$).

The picture is more complex than for the pregroup associated to the Leary-Stancu group. So far we don't know a characterisation of $P$ being a locality. An interesting candidate for $\Delta$ would be to take all the $O_p(G_i)$ for $i \in \{1, \cdots ,r\}$, together with their $S$-conjugate and the overgroups of all these in $S$. We then need to assume $S_i = N_S \big( O_p(G_i) \big)$ (otherwise it is easy to construct a triple $(g,s,g) \in \Db_\Delta$ but not in $\Db_P$). One can then prove that $\Db_P \subseteq \Db_\Delta$. However the converse inclusion does not hold in general: for example we need that any $R \leq S$ such that $O_p(G_i) \leq R \leq T_{i,g}$ for an $i \in \{ 1, \cdots , r\}$ and a $g \in H_i$ satisfies $N_S(R) \leq S_i$ (otherwise the word $(g^{-1},s,g)$ with $s \in N_S(R) \setminus S_i$ is in $\Db_\Delta$ via $({}^gR,R,R,R^g)$ but not in $\Db_P$). There are other necessary conditions coming from the fact that $\Db_P$ does not contain any words of the form $(aga',s,bhb')$ with $g$ in some $H_i$, $h$ in some $H_j$ and $i \neq j$. In any event, this seems to significantly reduce the possibilities for $P$ to be a locality.

\bibliography{biblio}{}
\bibliographystyle{plain}

\end{document}